%% file: PINN_Dis.tex
\documentclass[a4paper,
10pt
]{article}

\usepackage[english]{babel}

\usepackage[a4paper,left=3cm,right=3cm,top=2.5cm,bottom=2.5cm]{geometry}

\usepackage{amsfonts} 
\usepackage{hyperref}
\usepackage{bm}
\usepackage{authblk}
\usepackage{siunitx}
\usepackage{graphicx}
\usepackage[english]{babel}
\usepackage{amsmath}
\usepackage{xcolor}
\usepackage{subcaption}
\usepackage{float}
\usepackage{multirow}
{\color{blue}}

\usepackage{comment}
\usepackage{stackengine,graphicx}
\usepackage{tikz}
\usepackage{tabularx}
\usepackage[linesnumbered,ruled,vlined]{algorithm2e}

\title{Physics Informed Neural Network Framework for Unsteady Discretized Reduced Order System}
\author[1]{Rahul Halder\footnote{rhalder@sissa.it}}
\author[2]{Giovanni Stabile\footnote{giovanni.stabile@uniurb.it}}
\author[1]{Gianluigi Rozza\footnote{grozza@sissa.it}}

\affil[1]{Mathematics Area, mathLab, SISSA, via Bonomea 265, I-34136 Trieste, Italy}
\affil[2]{Department of Pure and Applied Sciences, Informatics and Mathematics Section, University of Urbino Carlo Bo, Piazza della Repubblica, 13, I-61029 Urbino, Italy}
\date{} 

\begin{document}
\maketitle
\newcommand{\GS}[1]{\noindent \textcolor{red}{GS: #1}}
\input{Abstract}
\input{Introduction}

\input{Governing_Equations}

\input{Physics_Informed_Neural_Netowork}
\input{Results}

\newpage
\input{Conclusion}
\input{Disclosure_statement}
\bibliographystyle{plain}

\bibliography{mybib}
\end{document}

%% file: Abstract.tex
\begin{center}
\bf ABSTRACT
\end{center}

This work addresses the development of a physics-informed neural network (PINN) with a loss term derived from a discretized time-dependent reduced-order system. In this work, first, the governing equations are discretized using a finite difference scheme (whereas, any other discretization technique can be adopted), then projected on a reduced or latent space using the Proper Orthogonal Decomposition (POD)-Galerkin approach and next, the residual arising from discretized reduced order equation is considered as an additional loss penalty term alongside the data-driven loss term using different variants of deep learning method such as Artificial neural network (ANN), Long Short-Term Memory based neural network (LSTM). 
The LSTM neural network has been proven to be very effective for time-dependent problems in a purely data-driven environment. The current work demonstrates the LSTM network's potential over ANN networks in physics-informed neural networks (PINN) as well. The major difficulties in coupling PINN with external forward solvers often arise from the inability to access the discretized forms of the governing equation directly through the PINN solver and also to include those forms in the computational graph of the network. This poses a significant challenge especially when a gradient-based optimization approach is considered in the neural network. Therefore, we propose an additional step in the PINN algorithm to overcome these difficulties. The proposed methods are applied to a pitch-plunge airfoil motion governed by rigid-body dynamics and a one-dimensional viscous Burgers' equation. The potential of using discretized governing equations instead of continuous form lies in the flexibility of input to the PINN. The current work also demonstrates the prediction capability of various discretized-physics-informed neural networks outside the domain where the data is available or governing equation-based residuals are minimized. 

\maketitle


%% file: Introduction.tex
\section{\bf{Introduction}}\label{sec:Intro}

Reduced order modelling \cite{rozza2022advanced,quarteroni2015reduced,benner2020model} for computational physics problems has gained
significant research focus in the past few decades. Usually, the computational effort of reduced-order model development involves two distinct stages - one is the offline stage where the high-fidelity simulation datasets are generated using high-performance computing (HPC) and the other one is the online stage which is performed using comparatively less powerful machines. Based on the learning methodology of the parametric space or unsteady dynamics,  the reduced order model can
be broadly classified into two categories, i.e. intrusive and non-intrusive reduced order model. In the intrusive approach, the governing equations are first projected on a reduced space, and the resulting equations are solved
numerically as a time-marching problem on the reduced space. Therefore, the intrusive schemes require significant modification of the solver whereas the
non-intrusive reduced order model, on the other hand, is entirely data-driven and a surrogate model is generated with a database from a high-fidelity numerical
or experimental dataset. Therefore, the second approach requires minimal changes in the computational solver.

The projection-based intrusive reduced order model often possesses numerical stability issues if the solution manifold is highly nonlinear, requiring a large dimension of the linear subspace used to approximate the solution manifold and increasing the computational cost of the
reduced order schemes. Among several intrusive reduced-order model approaches,
POD-Galerkin-based ROM \cite{STABILE2018273,KaratzasStabileNouveauScovazziRozza2019b,StabileHijaziMolaLorenziRozza2017,stabile2020efficient}, has been very effective in approximating the dynamics in a low dimensional manifold. A
survey of projection-based ROM for parameterized dynamical systems is carried
out by Benner et al.\cite{BennerROM}.  An additional effort is required in the
projection-based intrusive ROM approach to handle the nonlinear term requiring
several hyper-reduction methods \cite{BARRAULT2004667,BDEIM2010,carlberg2013gnat,farhat2015structure}. Since in an intrusive ROM framework, the governing equations are satisfied in the reduced or latent space, this type of ROM is expected to produce accurate prediction even at a sparse training dataset. 

There have been several progresses in non-intrusive
reduced order model in recent times \cite{benner2021system,demo1,demo2019complete,GadallaCianferraTezzeleStabileMolaRozza2020}. For ill-posed problems such as inverse, regression, and classification problems in diverse fields of computational mechanics, as outlined in Schmidhuber et al., \cite{schmidhuber2015deep}, several deep learning (DL) approaches offer effective solutions. Maulik and San \cite{maulik2019subgrid} and Duraisamy et
al.\cite{parish2016paradigm} proposed a turbulence closure model using machine learning. Among other notable applications of deep learning in fluid mechanics problems are the combined proper orthogonal decomposition (POD)
and LSTM network for incompressible flows with spectral proper orthogonal
decomposition (SPOD) of Wang et al. \cite{wang2018model} and the application of
a dimensionality reduction method with DL networks for learning feature dynamics
from noisy data sets in Lui and Wolf \cite{lui2019construction}. As mentioned,
non-intrusive approaches often require a large data set to gain significant
prediction accuracy. LSTM network has been applied to several fields of
computational physics problems such as aeroelastic applications and hydrodynamic
applications by Halder et al. \cite{AHalder2020,halder2023deep}, for the
application in atmospheric turbulence by Rahman et
al.\cite{rahman2019nonintrusive}. However, such approaches, being completely data-driven or black-box methods, require large training datasets computed in the expensive offline stage, and thereby they appear to be ineffective in those application areas where sparse datasets are available. Hence, in the non-intrusive framework, a grey-box approach where the governing dynamics or the parametric space can be learnt from the governing equations in addition to data-driven constraints can circumvent the requirement of a large training dataset. 

In the past century, early work on solving partial differential equations
(PDE) using neural networks starts with Dissanayake and Phan-Thien
\cite{dissanayake1994neural}. With recent advancements in research of neural networks, this
topic got renewed attention with the work of Raissi et al.\cite{CRAISSI2019686}.
Physics-informed neural networks can circumvent the problem of large training
dataset requirements by introducing additional physics-based loss terms. This
physics-based loss term in neural networks can be computed from the governing
equations. Different types of neural networks can be considered in the context of PINN, such as artificial neural networks (ANN) by Raissi et al. \cite{CRAISSI2019686}, Waheed et al.\cite{bin2021pinneik}, Tartakovsky et al. \cite{tartakovsky2020physics} and
Schiassi et al.\cite{schiassi2021extreme}. Similarly, Gao et al. \cite{gao2021phygeonet} and Fang \cite{fang2019physics} coupled convolution neural networks (CNN) with physics constraints. A review of several developments in PINN is summarized by Cuomo et
al.\cite{cuomo2022scientific}. However, there are limited works in Physics Informed Long Short-Term Neural Network (LSTM). Zhang et
al.\cite{cheng2021deep} used a tensor differentiator to determine the derivatives of state space variables to couple the LSTM network with physics information arising from the governing equations and applied to nonlinear-structural problems. PINN usually considers an automatic differential (AD) based approach \cite{baydin2018automatic} for the computation of the gradients. There have been few efforts to use derivatives
based on known numerical approaches such as finite difference, finite-element,
finite volume methods etc. to compute the physics-driven loss-term in the PINN
network. Ranade et al. \cite{RANADE2021113722} proposed a Discretization-net
where finite volume-based numerical residuals are considered in the loss term of
the physics-informed neural network. Similar work by Aulakh et al.
\cite{aulakh2022generalized} coupled physics-informed ANN network with the finite
volume-based discretization of the OpenFOAM solver. The PINN network often takes
more time than the forward solvers due to the large dimensionality of the network output. To reduce this effort the governing equation can be first reduced
using the conventional dimensionality reduction approach, such as Proper Orthogonal Decomposition (POD), and then can be coupled with the
loss equations as shown by Chen et al. \cite{CHEN2021110666} and
Hijazi et al. \cite{hijazi2023pod}.The full potential of the discretized physics-based neural network will be realized if one can couple the PINN solver and the external forward solver each independently formulated within distinct coding environments. In such a case, none of the previous works propose a seamless coupling of the external solver with the PINN environment without further modification of the forward solver. Discretized governing equation augmented
LSTM network is demonstrated for the first time by the author \cite{halder2023deep}, but the detailed comparison of prediction accuracy in the context of different neural networks and extension of the algorithm to a spatiotemporal problem will be demonstrated in the present work. 

The contribution of the current work can be summarized as follows:

\begin{itemize}
\item The current work proposes a novel approach where the existing forward solvers can be used in the PINN framework. Therefore, the solvers can be used for ill-posed and inverse problems. In the present work,
although we have coupled the PINN environment with a PyTorch-based external
numerical solver, the motivation of the current work is to bridge the PINN
environment with the existing C++-based open-source solvers for handling complex geometry and more practical applications as a future endeavour.

\item Current work proposes discretized governing equations based LSTM-PINN network for unsteady spatio-temporal problems of interest for the first time as per authors' knowledge.

\item Detail studies of prediction and reconstruction capabilities involving several neural networks such as ANN and LSTM and discretized governing equations are carried out.

\item The PINN network often suffers from the large-dimensionality of problems. Therefore, we propose to use the projection-based ROM framework used in \cite{STABILE2018273} in the PINN platform to accelerate the training time.

\item we propose an augmentation of the Physics-Informed Neural Network (PINN) algorithm, which incorporates an additional step to seamlessly integrate with any external solver. This additional step will eliminate the necessity for the incorporation of the discretized form of the governing equation within the computational graph of the neural network. 

\end{itemize}

The current paper is organized in the following way. 
\begin{itemize}
\item First, in \autoref{sec:GE}, the general formulation of the discretized governing equation is presented, followed by
the reduced-order representation of the discretized equation discussed. \item Next, in \autoref{sec:PINN}, the architecture of the ANN and LSTM network and their extension to a physics-informed neural network using additional loss penalty arising from the numerical residual of the reduced order discretized equations are discussed. In this section, we also discuss the software implementations \cite{demo2023extended} details of the work.

\item Furthermore, we propose an additional step of importing the physics-based residual and the Jacobian of the residual vector with respect to the output variable in the PINN solver if the external CFD solver is completely detached from the PINN environment to facilitate the inclusion of the discretized form of equation in the computational graph of the neural network. 

\item Finally, in \autoref{sec:results}, the proposed approaches are applied in the context of a 2-dof linear structural equation and a viscous Burgers' equation (full order and reduced order) to demonstrate the discretized physics-informed ANN and LSTM network. Next, in the context of the full-order Burgers' equation, we demonstrate our additional step of the PINN algorithm for seamless coupling of two separate environments - PINN solver and external forward solver. 

\end{itemize}

%% file: Governing_Equations.tex
\section{\bf{Governing Equations and Reduced Order System}}\label{sec:GE}
The governing equations associated with different computational mechanics' problems can be cast into the general form of a nonlinear parameterized dynamical system as shown in \autoref{eq:1} 
\begin{equation}\label{eq:1}
\begin{aligned}
&\frac{\partial {u}(t;\mu)}{\partial t} = \textbf{A} (\mu)  {u} (t;\mu) +
\textbf{F} ({u};\mu;t)  {u} (t;\mu) +
\textbf{B} (\mu) {f} (\mu), \ \
{u}(0;\mu) = {u}^{0}(\mu),
\end{aligned}
\end{equation}
where, $\mu \in {\Omega}_{\mu} \subset \mathbb{R}^{N_{\mu}}$ is a vector containing the
parameters associated with the dynamical system, $u:[0, T] \times
\mathbb{R}^{N_{\mu}} \rightarrow \mathbb{R}^{N} $ denotes the time-dependent state variable
function. $ {u}^0 : {\mathbb{R}}^{N_{\mu}} \rightarrow \mathbb{R}^{N} $ is the
initial state. $f:[0, T] \times \mathbb{R}^{N_{\mu}} \rightarrow {\mathbb{R}}^{N_i} $ denotes
the time-dependent input variable function independent of the state variable
function. $N$, $N_{i}$ and $N_{\mu}$ are the number of degrees of freedom of high-fidelity solution, dimensions of the input variable and the parameter space respectively. The constant real-valued operator ${\bf{A}}:{\mathbb{R}}^{N_{\mu}}
\rightarrow {\mathbb{R}}^{(N \times N)} $ contributes to the linear part of the
governing equation and the operator ${\bf{F}}: {\mathbb{R}}^{N_{\mu}} \rightarrow
{\mathbb{R}}^{(N \times N)} $ associated with the nonlinear part of the
governing equation is dependent on the state variable. The operator ${\bf{B}}:{\mathbb{R}}^{N_{\mu}} \rightarrow {\mathbb{R}}^{(N_i \times N_i)} $ associated
with the input variable can also be dependent on the
state variable based on the nature of the governing equation. Although, any time integration method can be
adopted to time-discretize the
dynamical system \autoref{eq:1}, we show a backward Euler time-discretization with an implicit approach as the numerical schemes for the demonstration purpose,  and the $\mathbf{B}$ matrix is considered here to be independent of the state variables.

\begin{equation}\label{eq:2}
\begin{aligned}
& (I_{N} - \Delta{t}^{(k)}A(\mu) - \Delta{t}^{(k)}F({u}^{(k)};\mu;t)) u^{(k)} = u^{(k-1)} + \Delta{t}^{(k-1)}B(\mu;t) f (\mu) , 
\end{aligned}
\end{equation}

whereas if we consider the explicit numerical discretization approach,

\begin{equation}\label{eq:3}
\begin{aligned}
& u^{(k)}  = u^{(k-1)} + (\Delta{t}^{(k-1)}A(\mu) + \Delta{t}^{(k-1)}F({u}^{(k-1)};\mu;t)) u^{(k-1)} + \Delta{t}^{(k-1)}B(\mu;t) f (\mu) , 
\end{aligned}
\end{equation}

where, $I_N$ denotes the identity matrix of size ${\mathbb{R}}^{(N \times N)}$,
$\Delta t$ is the time step and and $k$ denotes the $k^{th}$ time instant. The
following subsection will discuss the reduced-order formulation of the
discretized governing equations.

\subsection{Linear subspace solution representation}

The computational complexity of \autoref{eq:2} and \autoref{eq:3} is of order
$\mathcal{O}(N)$. Therefore, to reduce the computational effort of the full-order
system, the state variables can be represented as the linear combination of a
few basis vectors $\Phi$ and the projected equation can be solved in a
similar time-marching way as considered for the full-order model solution. We
will propose in the following \autoref{sec:PINN} how we can solve the set of
reduced-order equations in an alternative way in the context of physics-informed
neural network. To find the basis vectors from a set of solution vectors
associated with parameters $\mu$:

\begin{equation} \label{eq:4}
\begin{aligned}
& u (\mu) \approx {\Tilde{u}} (\mu) \equiv {\Phi} {\hat{u}}(\mu), 
\end{aligned}
\end{equation}

where, ${\hat{u}}(\mu): R^{N_{\mu}} \rightarrow R^{n}$ with $n << N$. The basis
vector sets $\Phi \in R^{N \times n}$ are defined as :

\begin{equation} \label{eq:5}
\begin{aligned}
& \Phi = \left[{\phi_1} \ \  ... \ \ {\phi_2}\ \ ...\ \ {\phi_{n}}\right].
\end{aligned}
\end{equation}
 To compute the basis vector sets, we collect the solution vectors $u$ at different time instants associated with the set of $m$ parameter values 
 $[{\mu_1,\mu_2,..\mu_{m}}]$. At parameter $\mu_p$, the snapshot vectors associated with total time instants of $N_t$ are
 defined as $U_p = [u^{(1)}(\mu_p), u^{(2)}(\mu_p),\ \ ....\ \ u^{(N_t)}(\mu_p)]
 \in R^{(N\times N_t)}$. Now, if we concatenate the snapshots
 corresponding to all the parameters, we obtain the total snapshot matrix $U
 \equiv [U_1 \ \ ... \ \ U_{N_\mu}]$, where $U \in R^{(N \times mN_t)}$. 


 Proper Orthogonal Decomposition (POD) \cite{STABILE2018273} can be used to
 compute the spatial basis vectors $\Phi$ by choosing the first $n$ columns of
 the left singular value matrix $W$ after the application of Singular Value
 Decomposition (SVD) on the snapshot matrix $U$ as shown in \autoref{eq:6}.

\begin{equation} \label{eq:6}
\begin{aligned}
& U = W \Sigma V^{T}, \\
&   =  \sum_{i=1}^{l} ({\sigma}_i w_i {v_i}^T),
\end{aligned}
\end{equation}

where $l = \min(N,mN_t)$ 
and $n < mN_t$. Here, $W \in \mathbb{R}^{N
\times l}$ and $V \in \mathbb{R}^{N_{\mu} N_t \times l}$ are set of orthogonal
vectors and $W \in \mathbb{R}^{l \times l}$ is a diagonal matrix where its
diagonal elements consist of the singular values. The POD basis ${\Phi}$
minimizes $\|U-{\Phi}{{\Phi}}^TU\|_F$ where $\|.\|_F$ denotes the Frobenius norm.
The reduced-order system of equations
can be written from the \autoref{eq:2} 

\begin{equation} \label{eq:7}
\begin{aligned}
({{\Phi}}^T I_{N}{{\Phi}} - \Delta{t}^{(k)}{{\Phi}}^T {\bf{A}}(\mu) {{\Phi}} - \Delta{t}^{(k)} {{\Phi}^T} {\bf{F}}({{\hat u}}^{(k)} ;\mu;t) {{\Phi}} ) {{\hat u}}^{(k)} = \\
{{\hat u}}^{(k-1)} +
\Delta{t}^{(k)}{{\Phi}}^T {\bf{B}}(\mu;t)  f (\mu),
\end{aligned}
\end{equation}

which can be written as follows:

\begin{equation}\label{eq:8}
\begin{aligned}
& (I_{n} - \Delta{t}^{(k)}A_R(\mu) - {{\Phi}^T} F({\Phi}{{\hat u}}^{(k)} ,\mu,t) {\Phi} {{\hat u}}^{(k)} = {{\hat u}}^{(k-1)} + \Delta{t}^{(k)}B(\mu;t) f (\mu) , 
\end{aligned}
\end{equation}

where, $A_R = {\Phi}^T A(\mu) {\Phi}$. Furthermore, the reduced order
system can be derived from the explicit formulation as shown in
\autoref{eq:3}. Similar to linear projection approaches such as Proper Orthogonal Decomposition (POD) as mentioned here, one can consider nonlinear manifold-based projection as shown by \cite{romor2023non} using the Convolutional Autoencoders approach.

\subsection{Test Cases} \label{sec:testcases}
As mentioned earlier, the current work considers two test cases: the
rigid body dynamical system for 2-dof pitch-plunge airfoil and the viscous
Burgers' equation to demonstrate the proposed methodologies. 

\textbf{(a) 2-dof pitch-plunge System} 
An implicit
second-order backward Euler time discretization approach is considered to
discretize the mass-spring system. The governing
equation of the mass-spring system can be written as follows: 

\begin{equation}\label{eq:9}
\begin{aligned}
\left[\begin{array}{cc}
1 & x_\alpha \\
x_\alpha & r_\alpha^2
\end{array}\right]\left(\begin{array}{l}
\ddot{h} \\
\ddot{\alpha}
\end{array}\right)+\left[\begin{array}{cc}
\left(\omega_h / \omega_\alpha\right)^2 & 0 \\
0 & r_\alpha^2
\end{array}\right]\left(\begin{array}{l}
h \\
\alpha
\end{array}\right)=\frac{V^{* 2}}{\pi}\left(\begin{array}{c}
-C l \\
-2 C m_{E A}
\end{array}\right).
\end{aligned}
\end{equation}

The terms $x_\alpha$, $r_\alpha$, $\omega_h$, $\omega_\alpha$ are structural
variables as mentioned in Halder et al.\cite{halder2022computational}. $h$ and
$\alpha$ are the plunge and pitch degrees of freedom, respectively, and $C_l$ and
$C_m$ are the lift and moment coefficient about the elastic axis
\cite{halder2022computational} typically used for aeroelastic applications. The
governing \autoref{eq:9} can be written as follows in matrix format:

\begin{equation}\label{eq:10}
\begin{aligned}
\mathbf{M} \ddot{X}+\mathbf{K} X = {F},
\end{aligned}
\end{equation}

where, $\mathbf{M}$ and $\mathbf{K}$ are the mass and stiffness matrix and the
variable $X$ is $[h \ \ \alpha]^T$ and $\bf{f}$ is force matrix arising from \autoref{eq:9} respectively. The residual computed from Eqns. \autoref{eq:9} and
\autoref{eq:10} can be written as follows:

\begin{equation}\label{eq:11}
\begin{aligned}
& r_\text{rigid-body} = \frac{3\left\{\begin{array}{c}
X \\
\dot{X}
\end{array}\right\}^{n+1}-4\left\{\begin{array}{c}
X \\
\dot{X}
\end{array}\right\}^n+\left\{\begin{array}{c}
X \\
\dot{X}
\end{array}\right\}^{n-1}}{2 d t} \\
& +\left[\begin{array}{cccc}
0 & 0 & -1 & 0 \\
0 & 0 & 0 & -1 \\
M^{-1} K & 0 & 0 & 0
\end{array}\right\}\left\{\begin{array}{c}
X \\
\dot{X} \\ 
\end{array}\right\}^{n+1} -\left\{\begin{array}{c}
0 \\
0 \\
M^{-1} F 
\end{array}\right\}^{n+1}.
\end{aligned}
\end{equation}


For the rigid-body dynamical systems, the number of dof is 2 and therefore, it does not require any further spatial reduction using \autoref{eq:4}. 

\textbf{(b) Viscous Burgers' Equation} 
An explicit discretization of the Burgers' equation with an upwind scheme for the convective term, the central difference scheme for the diffusion terms, and a forward difference scheme for the time derivative term is considered. For a given field $u(x,t)$ and kinematic viscosity $\nu$, the general form of the Burgers' equation can be written as $\frac{\partial u}{\partial t}+u
\frac{\partial u}{\partial x}=\nu \frac{\partial^2 u}{\partial x^2},$. The residual arising from the governing equation after discretization is as \autoref{eq:12}
\begin{equation}\label{eq:12}
\begin{aligned}
r_\text{burgers}=\frac{u_i^{n+1}-u_i^n}{d t}-u_i^n \frac{\left(u_{i+1}^n\right)-\left(u_i^n\right)}{d x}+\nu \frac{u_{i+1}^n-2 u_i^n+u_{i-1}^n}{(d x)^2},
\end{aligned}
\end{equation}

The residual, $r_{burgers}$ can be expressed using a reduced order variable
$\hat{u}$:

\begin{equation}\label{eq:13}
\begin{aligned}
r_\text{burgers,reduced}=\frac{{\hat{u}}_i^{n+1}-{\hat{u}}_i^n}{d t}- {\Phi}^T {\Phi}\hat{u}_i^n \frac{\left(\Phi \hat{u}_{i+1}^n\right)-\left(\Phi \hat{u}_i^n\right)}{d x}+ \nu \frac{{\hat{u}}_{i+1}^n-2 {\hat{u}}_i^n+{\hat{u}}_{i-1}^n}{(d x)^2},
\end{aligned}
\end{equation}
where $dt$ and $dx$ are the time-step and spatial discretization steps suitable
for the numerical discretization approach in the current test case, and $\nu$ is
the viscosity coefficient. 

\subsection{Hyper Reduction} \label{sec:hyperreduction}
As demonstrated by Chaturantabut et al.\cite{BDEIM2010}, the complexity of the nonlinear term 
is not reduced in \autoref{eq:8}. The nonlinear part of the reduced order equation can be written as :
\begin{equation}\label{eq:14}
\begin{aligned}
\tilde{\mathbf{N}}(\hat{u}):=\underbrace{{\Phi}^T}_{n \times N} \underbrace{\mathbf{F_{nl}}\left(\mathbf{\Phi} \hat{u}(t,\mu)\right)}_{N \times 1} .
\end{aligned}
\end{equation}

where the nonlinear term ${\bf{F_{nl}}} = F({\Phi}{{\hat u}} ,\mu,t) {\Phi} {{\hat u}}$. The nonlinear term has a computational complexity of $(\mathcal{O}(\alpha(N)+4
Nn))$, where $\alpha$ is some function of
$N$. Empirical Interpolation Method (EIM) as proposed by Barrault et al.\cite{BARRAULT2004667}, is a popular hyper-reduction approach to reduce the computational complexity of the nonlinear term. Hence in the current work, a variant of EIM approach,
Discrete Empirical Interpolation Method (DEIM) as mentioned by Chaturantabut and
Sorensen \cite{BDEIM2010} is adapted to handle the nonlinear term. However, any
other hyper-reduction approach such as the GNAT method by Carlberg et
al.\cite{carlberg2013gnat}, the ECSW method by Farhat et
al.\cite{farhat2015structure} can also be considered. In the hyper-reduction
approaches, the nonlinear term, $\mathbf{F_{nl}} (t,\mu)$ can be reconstructed from the
information at a few collocation points in the computational domain. First, the
nonlinear term $\mathbf{F_{nl}}(t,\mu)$ can first be projected on a subspace which
approximates the space generated by the nonlinear term and spanned by the basis
of dimension $m_h << N$. The basis vector sets ${\Phi^h} = \left\{{\phi^h}_1,
\ldots, {\phi^h}_{m_h}\right\}$. 
The nonlinear term then $\bf{F_{nl}}$ can be
approximated as ${\bf{F_{nl}}}(t, \mu) = {{\Phi^h}} \mathbf{c}(t,\mu)$, where the
coefficient vector $\mathbf{c}(t,\mu)$ can be computed from an index matrix or
mask matrix $\mathbf{P}$ used for indexing sensor locations in the computational
domain. The index matrix is defined as $\mathbf{P}=\left[e_{\rho 1}, \ldots
\ldots , e_{\rho l}\right] \in \mathbb{R}^{N \times m}$ where,
$\mathbf{e}_{\wp_i}=[0, \ldots, 0,1,0, \ldots, 0]^T \in \mathbb{R}^N$. $1$ is
placed at the $i^{th}$ entry of the $\mathbf{e}_{\wp_i}$ vector. $\mathbf{P}$
can be computed using the following relationship in \autoref{eq:15}:

\begin{equation}\label{eq:15}
\begin{aligned}
\mathbf{P}^T {\mathbf{F_{nl}}}(\hat{u}, t, \mu)={\left(\mathbf{P}^T{\Phi^h}\right)} \mathbf{c}(t,\mu),
\end{aligned}
\end{equation}

and, the nonlinear term $\mathbf{F_{nl}}(t, \mu)$ can be approximated as follows: 

\begin{equation}\label{eq:16}
\begin{aligned}
{\mathbf{F_{nl}}}(t, \mu)={{\Phi^h}}\left(\mathbf{P}^T {\Phi^h}\right)^{-1} \mathbf{P}^T \mathbf{F_{nl}}(\hat{u}, t, \mu),
\end{aligned}
\end{equation}

$\mathbf{P}^T \mathbf{F_{nl}}(t, \mu)$ indicates the location of the computational
domain where the nonlinear term needs to be computed. Therefore, to reduce the
computational complexity, the final nonlinear term to be used in the reduced
form of the discretized \autoref{eq:13} can be written as follows:
\begin{equation}\label{eq:17}
\begin{aligned}
\tilde{\mathbf{N}}(\hat{u})= {\Phi}^T {\Phi^h}\left(\mathbf{P}^T {\Phi^h}\right)^{-1} \mathbf{P}^T \mathbf{F_{nl}}(\hat{u}, t, \mu).
\end{aligned}
\end{equation}

%% file: Physics_Informed_Neural_Netowork.tex
\section{\bf{PINN for Discretized Reduced Order System}}\label{sec:PINN}

In this section, physics informed neural network using the discretized governing equation will be discussed. First, two different types of neural network
architecture such as ANN and LSTM will be briefly introduced in 
\autoref{subsec:NN} followed by the integration with the discretized governing equation in the framework of the physics-informed neural network
\cite{CRAISSI2019686} which will be discussed in \autoref{subsec:ANN-PINN}
and \autoref{subsec:LSTM-PINN}.  

\subsection{Neural Network Architecture}\label{subsec:NN}

ANN formulates a nonlinear functional relationship between the input and output.
In \autoref{eq:18}, $\textbf{f}$ can be considered as input whereas the variable
$\textbf{u}$ can be the output of the neural network and input-output mapping
for a single-layer, ANN network can be written as follows:

\begin{equation}\label{eq:18}
\begin{aligned}
\textbf{u} = f_{\text{act}}(W \textbf{f} + b),
\end{aligned}
\end{equation}

where $f_{act}$ is the nonlinear activation function. 
the weight matrix and
$b$ is the bias vector. In \autoref{eq:18}, the input, $\textbf{f}$ and output $\textbf{u}$ can also
correspond to an unsteady time series problem such as the mass-spring system and
unsteady viscous Burgers' equation. The Eqn.\autoref{eq:18} can be rewritten as \autoref{eq:19}:  

\begin{equation}\label{eq:19}
\begin{aligned}
p_o^{n_t} = f_{act}(W p_i^{n_t} + b),
\end{aligned}
\end{equation}
where ${p_i}^{n_t}$ corresponds to the input of the unsteady problem considered
which is $C_l$ and $C_m$ at $n_t^{th}$ time instant for the pitch-plunge airfoil and
$t$ for the Burgers' equation. Furthermore, $p_o^{n_t}$ corresponds to the output
of the unsteady problem considered which is $h$ and $\alpha$ at $n_t^{th}$ time
instant for the pitch-plunge airfoil and $u$ at all the spatial locations for
the Burgers' equation. However, in Recurrent Neural Networks (RNN), the output,
$p_o^{n_t}$ corresponds to the time history of the input instead of a single time
instant value as shown in \autoref{eq:20}. The other matrix $U_R$ is not present
in the architecture of ANN which is associated with the input from the
previous time steps values and introduces the memory effect in the network
architecture. In different variants of RNN, $h_i^{(n_t)}$ forms a different
functional relationship with previous time instant values of $p_i$. For a simple
RNN network, $h_i^{(n_t)} = p_i^{n_t-1}$   

\begin{equation}\label{eq:20}
\begin{aligned}
p_o^{n_t}=f_\text{act}\left(Wp_i^{n_t}+U_{R} h_i^{(n_t)}+b\right).
\end{aligned}
\end{equation}

\begin{figure}[ht] \label{inpoutLSTM}
\centering
\begin{subfigure}[b]{0.90\textwidth}
\centering
\subfloat[ Input Tensor of LSTM
]{\label{fig:Fig1a}\includegraphics[width=.60\linewidth]{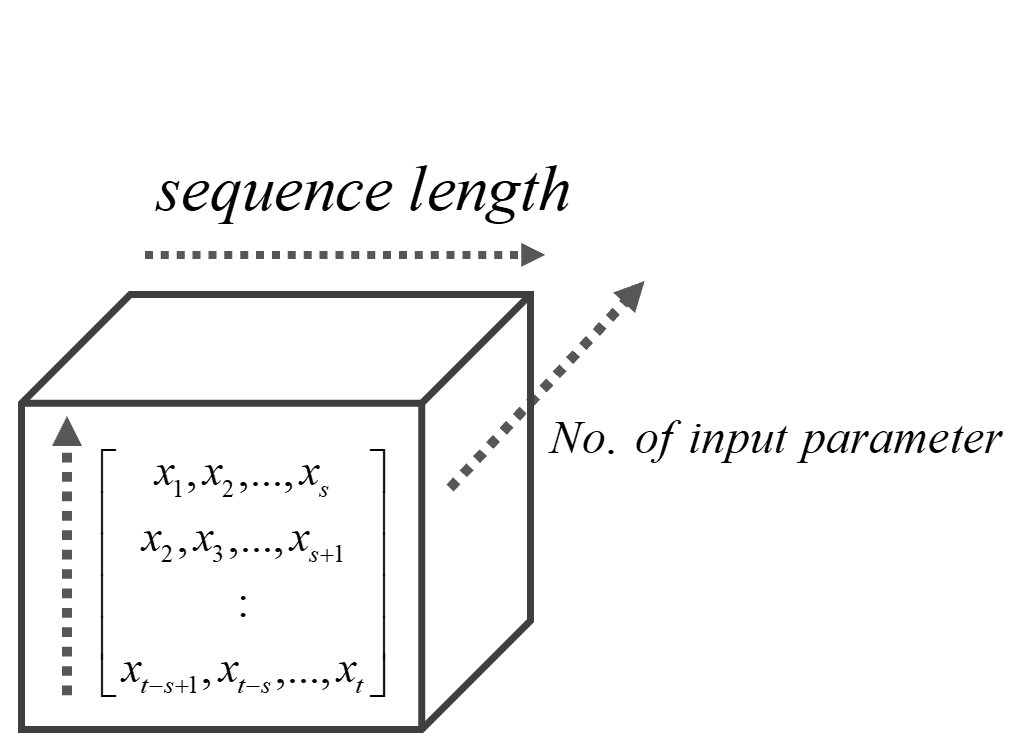}}
\subfloat[Output of
LSTM]{\label{fig:Fig1b}\includegraphics[width=.35\linewidth]{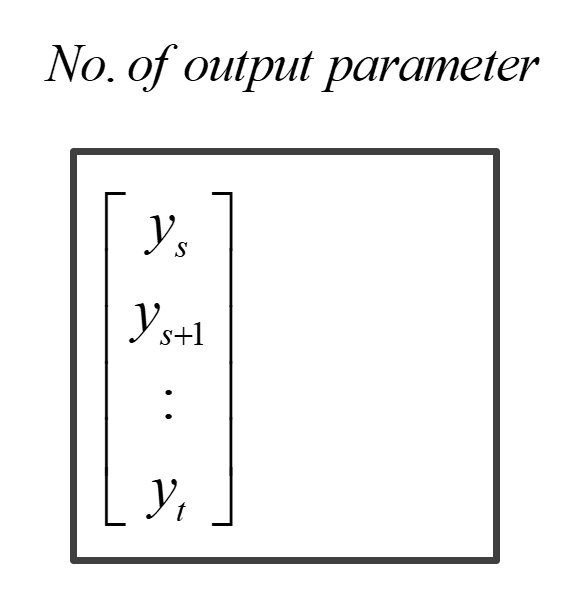}}
\end{subfigure}
\medskip
\caption{Input and output architecture of LSTM.}
\label{fig:inpoutLSTM}
\end{figure}

The present work considers a variant of the RNN network called a Long Short Term
Memory network (LSTM). A brief introduction to an LSTM cell and the associated
network's input and output architecture will be briefly introduced now, which
will be essential to understand the coupling between the LSTM network with the
discretized PDE-based loss term, a topic of discussion in the following
subsection. However, the advantages of the LSTM network over a simple RNN
network, and the detailed architecture of the LSTM network can be found in
Halder et al. \cite{AHalder2020,halder2023deep}. As mentioned, RNN-type
neural networks require input values associated with a sequence of time
instants. The output corresponds to the last time instant in the input time
sequence. For example, for an input time-series, $[x_1, x_2,x_3 .... x_t]$, corresponding output, $[y_1, y_2,y_3 .... y_t]$ and sequence length of $s$ associated with the input, in the LSTM network, $[x_1, x_2,...x_s]$ corresponds to output
$y_s$. The input of the LSTM forms a tensor structure as shown in
\autoref{fig:Fig1a} where the number of input variables attributed to the third
dimension. For example, for a pitch-plunge airfoil system, the number of input
variables is 2 which are $C_l$ and $C_m$ and for the Burgers' equation, time, $t$
being the only input variable, its dimension is $1$.
Conversely, the output structure forms a matrix as shown in \autoref{fig:Fig1b}
where the second dimension corresponds to the number of output variables. For a
pitch-plunge airfoil system, the number of output variables is 2 which are $h$
and $\alpha$, whereas, for the Burgers' equation the number of output variables
is the number of spatial locations.

However, this input-output structure of both ANN and LSTM networks is adequate
to understand the discretized reduced order system enhanced physics-informed
neural network termed ANN-DisPINN and LSTM-DisPINN respectively which we will discuss in the following 
\autoref{subsec:ANN-PINN} and \autoref{subsec:LSTM-PINN} .

\subsection{ANN-DisPINN}\label{subsec:ANN-PINN} 
\begin{figure}[t]
    \centering
    \includegraphics[width=100mm,scale=0.5]{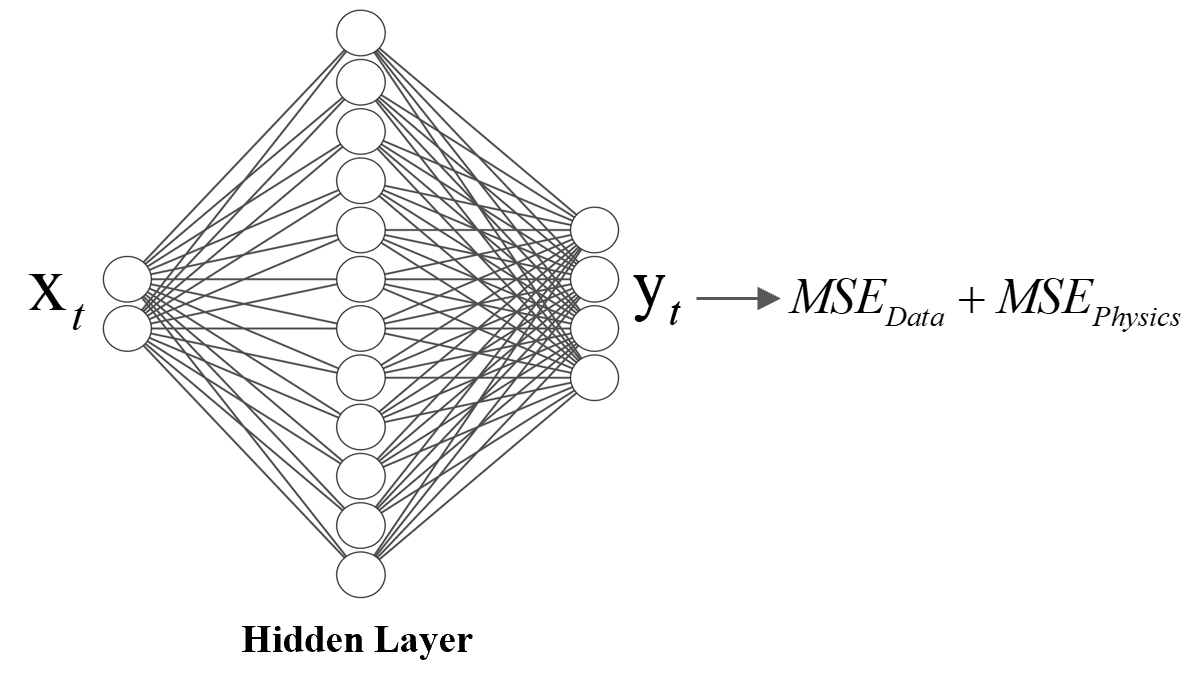}
        \caption{Architecture of ANN-DisPINN.}
        \label{fig:Fig2}
\end{figure}

In this section,  the discretized
full-order or reduced-order governing equation is employed in the
physics-informed neural network framework. The input and output structure of a
simple ANN and ANN-Network is straightforward. For an unsteady problem, $[x_1,
x_2,x_3,\dots,x_t]$ is the input associated with different time instants and
$[y_1, y_2,y_3,\dots,y_t]$ are the output values corresponding to the inputs.
Hence, a one-to-one correspondence exists between the input and output
formulations associated with each time instant. However, in ANN-DisPINN proposed
here in \autoref{fig:Fig2}, an additional loss penalty is introduced in addition to
the data-driven loss part, which is coming from the discretized governing
equation as shown in \autoref{eq:11} and \autoref{eq:13}. The data-driven and
physics-driven loss functions associated with the pitch-plunge airfoil and
viscous Burgers' equation are as follows:

\begin{equation}\label{eq:21}
\begin{aligned}
MSE_\text{Data} = || y_\text{pred} - y_\text{actual}||_{L2},\\
MSE_\text{burgers} =  || r_\text{burgers} ||_{L2}, \\
MSE_\text{rigid-body} =  || r_\text{rigid-body} ||_{L2},
\end{aligned}
\end{equation}

where, $y_\text{pred}$ is the predicted outcome from \autoref{eq:18} and $y_\text{actual}$
is the benchmark numerical dataset available. The loss terms arising from the
data-driven term $MSE_\text{Data}$ and physics-driven term $MSE_\text{burgers}$
and $MSE_\text{rigid-body}$ need to be minimized to obtain the weight matrix
$W$ and bias vector $b$ as mentioned in \autoref{eq:18} from this
non-convex optimization problem statement which is inherent to the neural network
algorithm. The 2-dof rigid body dynamical system does not require any spatial
reduction whereas, the Burgers' equation is first projected on a reduced basis
followed by residual computation from the reduced discretized equation which
further serves as a loss penalty in the physics-driven loss term of ANN network
as $MSE_\text{burgers-reduced} =  || r_\text{burgers,reduced} ||_{L2}$ 

\subsection{LSTM-DisPINN}\label{subsec:LSTM-PINN}

\begin{figure}[ht] 
\centering
\begin{subfigure}[b]{0.90\textwidth}
\centering
\subfloat[ A single LSTM Cell]{\label{fig:Fig3a}\includegraphics[width=.35\linewidth]{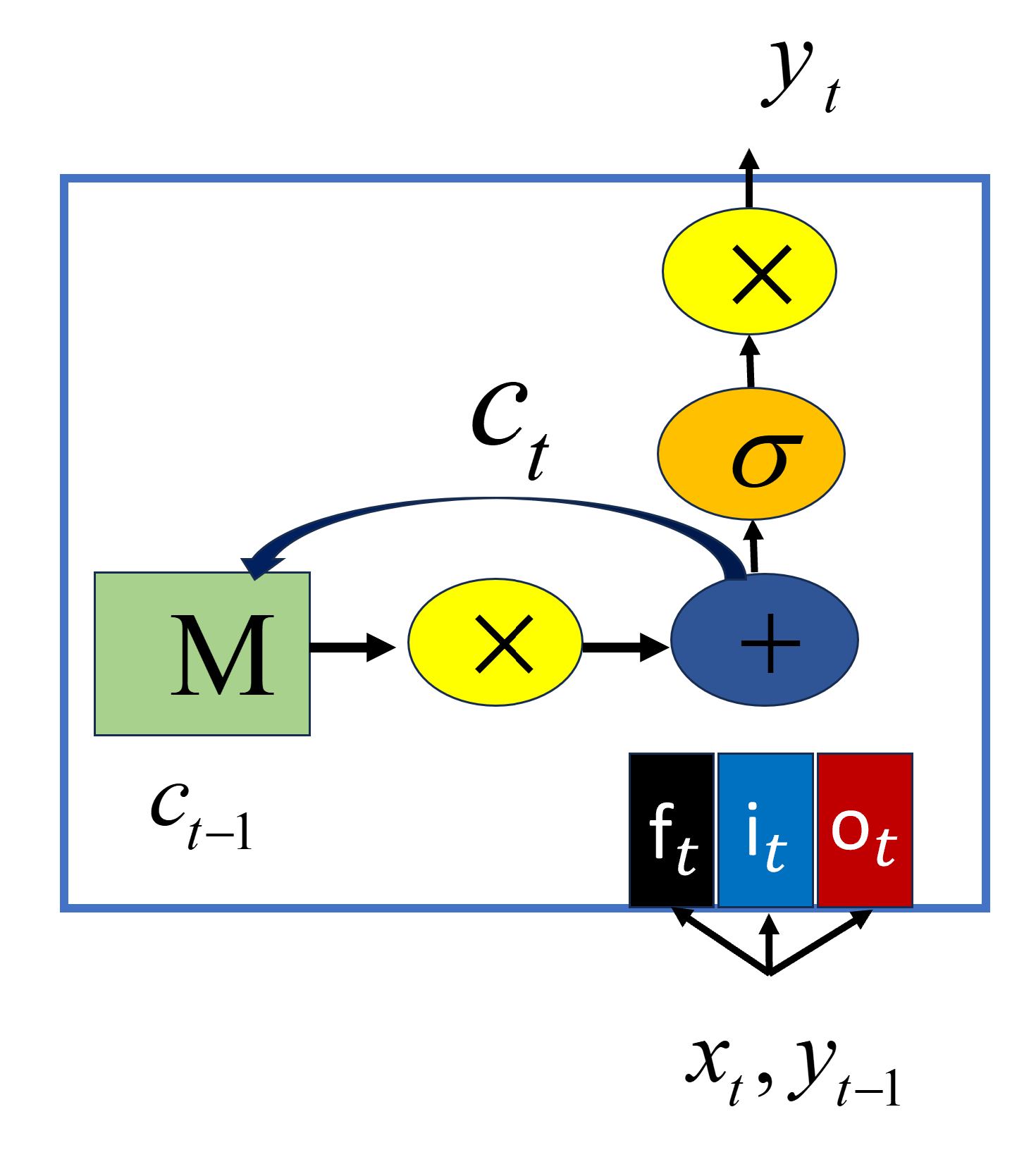}}
\subfloat[LSTM-PINN
]{\label{fig:Fig3b}\includegraphics[width=.65\linewidth]{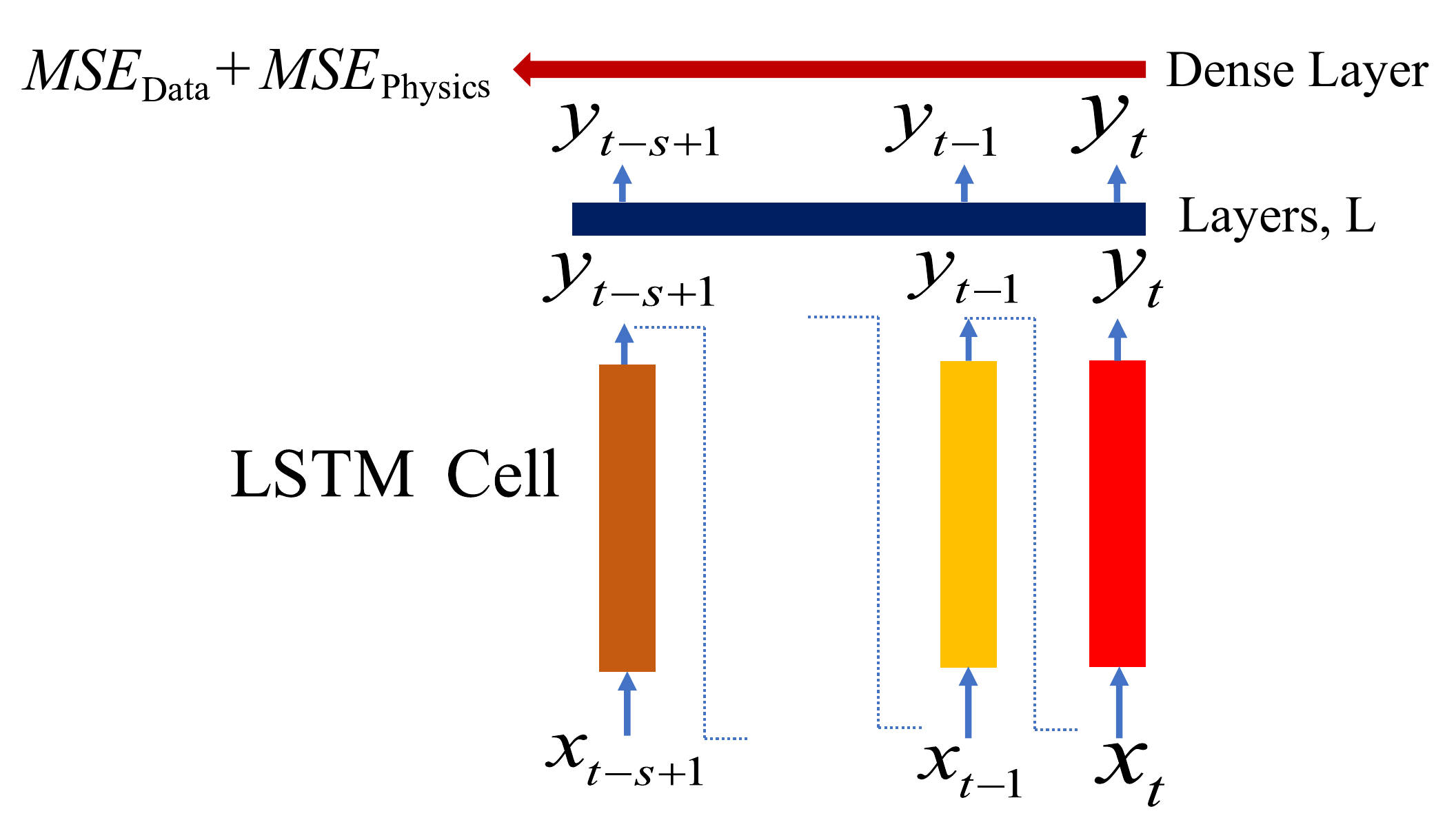}}
\end{subfigure}
\medskip
\caption{Architecture of a single LSTM-Cell and Discretized PDE-based LSTM-DisPINN.}\label{fig:Fig3}
\end{figure}

As shown in \autoref{fig:Fig1a}, the input to the network is
first re-formulated with a sequence length and then it is passed through the
LSTM architecture as shown in \autoref{fig:Fig3a} and \autoref{fig:Fig3b}  to
obtain the $y_\text{pred}$ from the nonlinear functional relationship between the input
and output. It is clear from \autoref{fig:Fig1b}, we only obtain predicted
output from $s^{th}$ time step. The predicted $y$ from the LSTM  as shown in \autoref{fig:Fig3b} is further used to minimize the
loss terms $MSE_\text{Data}$, $MSE_\text{physics, burgers}$ and $MSE_\text{physics,
rigid-body}$ to obtain the weight matrices and bias vectors associated with the
LSTM network. Similarly, for the reduced order discretized system, $MSE_\text{physics,
burgers-reduced}$ should be considered in the physics-driven loss term.  \\

The advantages of using the discretized system in the loss-term of the discretized-physics-informed neural networks are summarised below:
\begin{itemize}

  \item The automatic differentiation (AD) based gradient computation \cite{baydin2018automatic} in
  conventional PINN framework requires the spatial coordinates (i.e., Cartesian
  coordinates) and time coordinates as the input to the neural network where the current approach provides a flexible choice of input.  
  \item If a sparse numerical dataset is available to introduce in the
  data-drive loss term of the PINN network, they should comply with the
  numerical derivatives used in the physics-driven loss function, which have
  been used to generate those numerical datasets. However, in the conventional
  PINN framework, the numerical data often does not comply with the AD-based
  derivative resulting in a less accurate prediction.
\end{itemize}

\begin{figure}[t]
    \centering
    \includegraphics[width=1\linewidth]{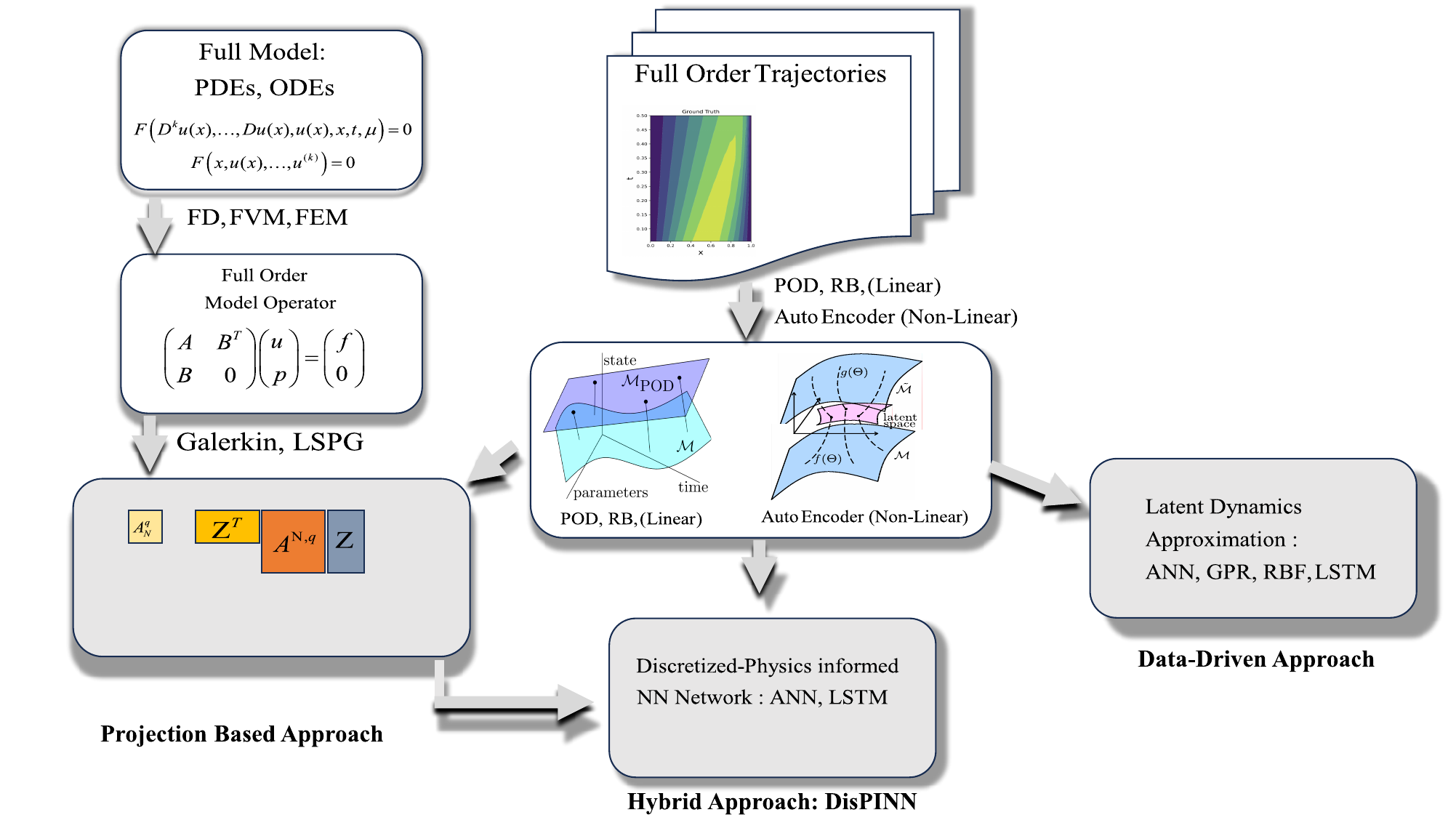}
        \caption{General framework for different reduced order models (i.e., projection-based approach, data-driven approach) and DisPINN, a bridge between these two aspects.}
        \label{fig:DisPINN}
\end{figure}

The discretized-physics-based neural network(DisPINN) can be put in the context of a projection-based reduced order model and data-driven approach as shown in \autoref{fig:DisPINN}. The implementation of the ANN-DisPINN and LSTM-DisPINN is carried out in the PyTorch-based PINN platform, PINA \cite{coscia2023physics}. 

\subsection{External Solver Coupling} \label{sec: external}

Here we discuss the coupling procedure of PINN with the external solver detached from the computational graph in the neural network. In a conventional neural network framework, the derivative of the physics-based loss term $(\text{L}_\text{eqn})$ and data-driven loss-term $(\text{L}_\text{Data})$ concerning the neural network parameters, $p = [W,b]$ is computed from the computational graph using backpropagation algorithm as mentioned by Amari et al.\cite{amari1993backpropagation}. The data-driven loss term is : 
\begin{equation}\label{eq:22}
\begin{aligned}
\mathrm{L}_{\text {Data }}=\frac{1}{N_{\text {Data }}} \sum_{i=1}^{N_{\text {Data }}}\left(y_{\text {pred }}-y_{\text {actual }}\right)^2
\end{aligned}
\end{equation}

Whereas the physics-driven loss term is : 

\begin{equation}\label{eq:23}
\begin{aligned}
\mathrm{L}_{\text{eqn}}=\frac{1}{N_{\text{eqn}}} \sum_{i=1}^{N_{\text{eqn}}} R^2
\end{aligned}
\end{equation}

where, the physics-based residual, $R$ can arise from full-order or reduced-order equations, such as $r_\text{burgers, reduced}$ and $r_\text{burgers}$. $N_{\text {Data}}$ and $N_{\text {eqn}}$ is the number of points available for the computation of $(\text{L}_\text{Data})$ and residual, $R$ respectively. The derivative of ${\text{L}}_\text{eqn}$ concerning the parameters can be computed as follows:

\begin{equation}\label{eq:24}
\begin{aligned}
\partial \mathrm{L}_{\text {eqn }} / \partial \mathrm{p}=\left(\frac{1}{N_{\text {eqn }}} \sum_1^{N_{\text {eqn}}} 2 R \frac{\partial R}{\partial p}\right)=\left(\frac{1}{N_{\text {eqn }}} \sum_1^{N_{\text {eqn }}} 2 R \frac{\partial R}{\partial y_{pred}} \frac{\partial y_{pred}}{\partial p}\right)
\end{aligned}
\end{equation}

Now, if the physics-based residual, $R$ is computed from an external solver, it is very difficult to access the discretized form directly from the PINN solver, therefore we can not compute $\partial \mathrm{L}_{\text {eqn }} / \partial \mathrm{p}$ from the computational graph using back-propagation. However if both the $R$ and $\frac{\partial R}{\partial U}$ are imported from the external solver, the modified physics-based loss therm, $\mathrm{L}_{D i s}$  that can be passed through the backpropagation is written as follows:
\begin{equation}\label{eq:LDIS}
\begin{aligned}
\mathrm{L}_{D i s}=\left(\frac{1}{N_{\text {eqn }}} \sum_1^{N_{\text {eqn }}}\left[2 R \frac{\partial R}{\partial y_{pred}}\right]_{\operatorname{detached}} y_{pred}\right)
\end{aligned}
\end{equation}

Since,  both the $R$ and $\frac{\partial R}{\partial U}$ are imported from the external solver, therefore detached from the computational graph, the derivative concerning the parameter is only computed on $y_{pred}$. Hence the loss term term which needs to be passed through the backpropagation is ($\text{L}= \mathrm{L}_{\text {Dis }}+\mathrm{L}_{\text {Data }}$). Finally, a stochastic gradient descent approach is considered in the current work as mentioned in \cite{amari1993backpropagation} to update the weight matrices and bias vectors based on $\partial \mathrm{L} / \partial \mathrm{p}$ following the relationship:

\begin{equation}\label{eq:pupdate}
\begin{aligned}
p_{new}= p_{old}- \alpha \frac{\partial \text{L}}{\partial p}
\end{aligned}
\end{equation}

where, $p_{new}$ and $p_{old}$ are the updated and old parameter values and $\alpha$ is the learning rate. If the residual $R$ is computed from the external solver, Jacobian, $J = \frac{\partial R}{\partial y_{pred}}$ needs to be updated at every epoch. However, $\frac{\partial R}{\partial y_{pred}}$ is a sparse matrix, and the overall computational expense of PINN solver can increase several folds if the $J$ term is updated every epoch and it contains residual from the entire computational domain. Therefore the computational expenditure can be reduced by following steps. 

\begin{itemize}
  \item Jacobian term, $J$ is updated only at a certain epoch interval.
  \item Residual, $R$ is computed at collocation points in the spatial domain which is an outcome of the reduced order method as mentioned in \autoref{sec:GE}
\end{itemize}

In the current work, the $J$ term is computed using finite-difference (FD), and no additional modification is required in the external solver, demonstrating the potential of our approach for seamless coupling of PINN with any other external forward solver. \autoref{alg:algo} summarize the proposed DisPINN in a pseudo-code. As indicated earlier, the primary inputs to the DisPINN algorithm are input to a prescribed neural network (NN), $\textbf{u}$, the number of points for the computation of physics-based and data-driven loss terms as $N_\text{eqn}$ and $N_\text{data}$. The iteration number in the optimization process is termed as $epoch$. $k$ is the $epoch$ interval at which the Jacobian, $J$ is updated. The updated $J$ after the prescribed interval is termed as $J_{new}$.   

\SetKwComment{Comment}{/* }{ */}

\begin{algorithm}[hbt!] 
\caption{DisPINN pseudo-code for coupling external detached forward solver}\label{alg:algo}
\KwData{$y_\text{actual}, \textbf{u}, N_\text{data}, N_\text{eqn}$}

$[{W},{b}] \gets \textbf{INIT} ([{W},{b}])$ \Comment*[r]{Initialize Weight and bias}
\While{$\text{epoch} \leq \text{Total Epoch No.}$}{
  
  $y_\text{pred} \gets \text{NN}({W},{b},f_\text{act}, \textbf{u})$ \Comment*[r]{prediction from chosen network $\text{NN}$}
  $\text{Ext. Solver} \gets y_\text{pred}$ \Comment*[r]{pass to external solver}
  $ R \gets \text{Ext. Solver}$  \Comment*[r]{pass to PINN solver from external solver}
  $\mathrm{L}_{e q n} \gets \frac{1}{N_{\text {eqn }}}\sum_{i=1}^{N_{e q n}} R^2$ \Comment*[r]{compute physics-based loss term}
  $\mathrm{L}_{\text {Data }} \gets \frac{1}{N_{\text {data }}}\sum_{i=1}^{N_{\text {data }}}\left(y_{\text {pred }}-y_{\text {actual }}\right)^2$ \Comment*[r]{compute data-based loss term}

  \eIf{$\text{rem}(epoch, k) = 0$}{  
    $\frac{\partial R}{\partial y_\text{pred}} \gets \text{Ext. Solver}$ \Comment*[r]{compute Jacobian from external solver}
    $J_\text{new} \gets \frac{\partial R}{\partial y_\text{pred}}$\; 
    $\mathrm{L}_{Dis} \gets (\frac{1}{N_{\text {eqn }}} \sum_1^{N_{\text {eqn }}}2 R J_\text{new} y_\text{pred})$  \Comment*[r]{compute additional physics based loss term for backpropagation}
    }
    {$J \gets J_\text{new}$ \;
    $\mathrm{L}_{D i s} \gets (\frac{1}{N_{\text {eqn }}} \sum_1^{N_{\text {eqn }}}2 R J y_\text{pred})$ \;
  }

  $L \gets (\mathrm{L}_{D i s} + \mathrm{L}_{\text {Data }})$ \Comment*[r]{compute total loss term for backpropagation}

  $\Delta {W} \leftarrow-\alpha \mathrm{G}_{\mathrm{ADAM}}\left(\nabla_{{W}} L\right), \Delta {b} \leftarrow-\alpha \mathrm{G}_{\mathrm{ADAM}}\left(\nabla_{{b}} L\right)$ \; \Comment*[r]{compute weight matrices and bias vector update}
  ${W} \leftarrow {W}+\Delta {W}, {b} \leftarrow {b}+\Delta {b}$ \;
}
\end{algorithm}

%% file: Results.tex
\section{\bf{Numerical results}}\label{sec:results}

This section considers two unsteady test cases: an unsteady pitch-plunge airfoil motion under a sinusoidal excitation and an unsteady viscous Burgers' equation. We first discuss the details of the high-fidelity models, followed by the surrogate models, ANN-DisPINN and LSTM-DisPINN corresponding to both test cases in \autoref{subsec:MS} and \autoref{subsec:Burgers_Full}. Finally, the PINN network based on the reduced order discretized equation is demonstrated in the context of the viscous Burgers' equation in \autoref{subsec:Burgers_reduced}. Since, in our current work, both the external solver and the PINN solver are developed in the PyTorch environment, it is not difficult to access the discretized form of the governing equation in the external solver from the PINN and include those forms in the computational graph for the computation of derivatives required for gradient descent optimization. However, in \autoref{subsec:new_loss}, the output residual and the Jacobian are detached from the computational graph before they are passed to the PINN solver to implement \autoref{alg:algo} and application is demonstrated on the full-order Burgers' equation.  

\subsection{Mass-Spring System} \label{subsec:MS}

\begin{figure}[ht]
\centering
\begin{subfigure}[b]{0.80\textwidth}
\centering
\subfloat[ Data Point =
1, AD-PINN]{\label{fig:loss_AD_1}\includegraphics[width=.50\linewidth]{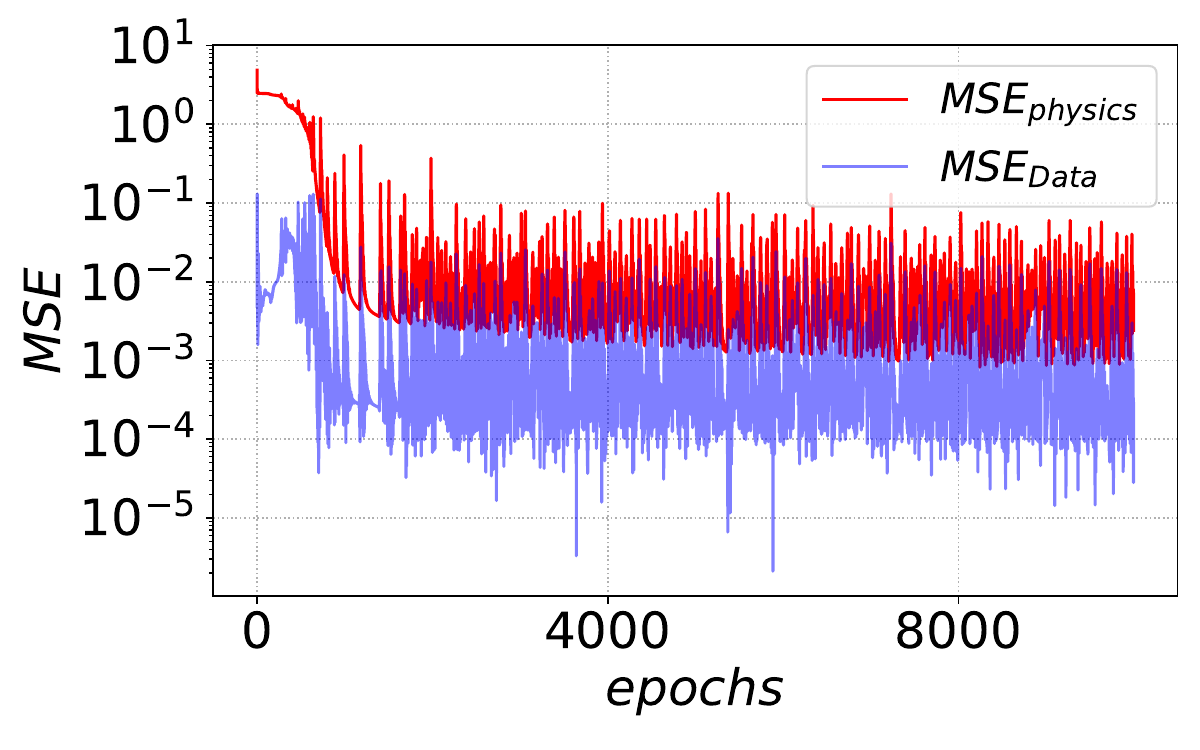}}
\subfloat[Data Point =
1000, AD-PINN
]{\label{fig:loss_AD_1000}\includegraphics[width=.50\linewidth]{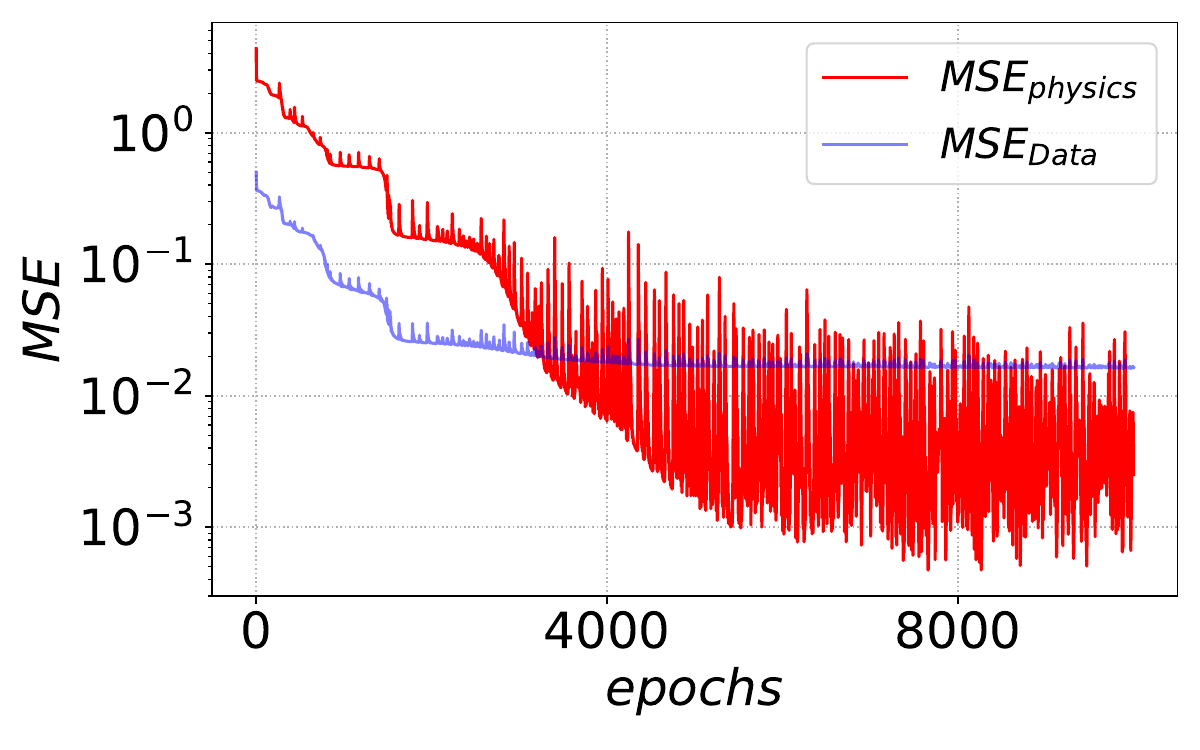}}
\end{subfigure}
\begin{subfigure}[b]{0.80\textwidth}
\centering
\subfloat[ Data Point =
1, ANN-DisPINN]{\label{fig:loss_ANN_1}\includegraphics[width=.50\linewidth]{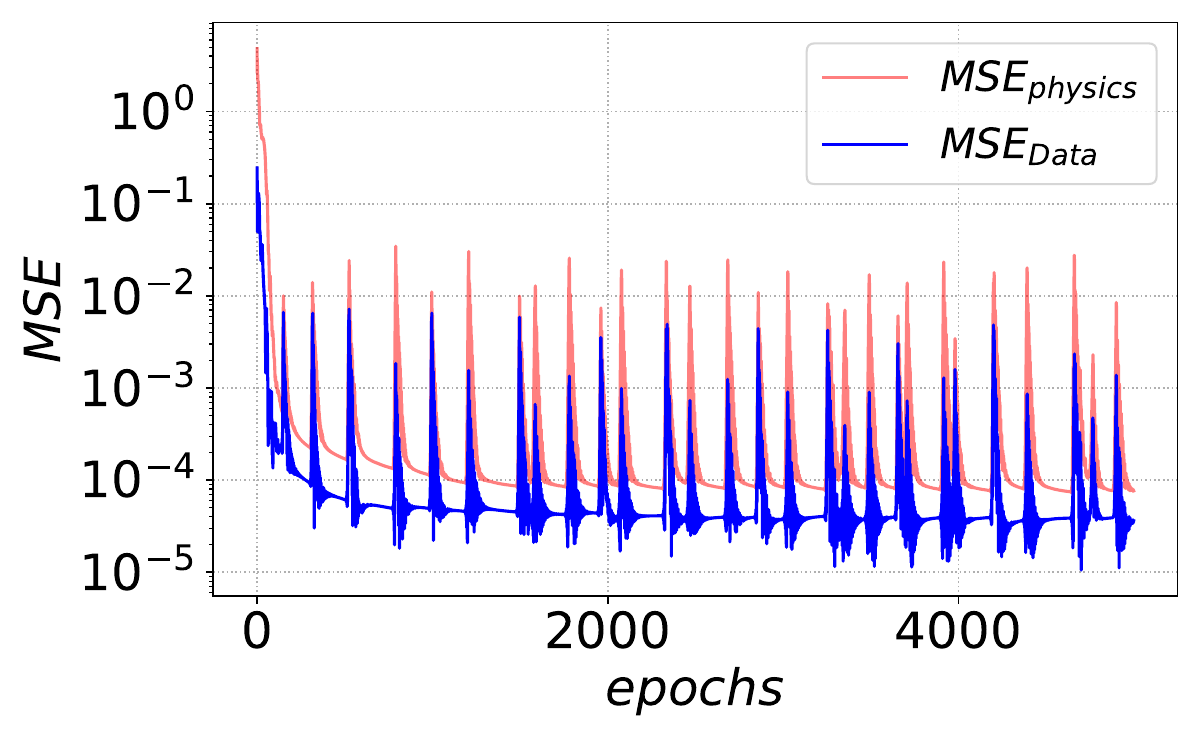}}
\subfloat[ Data Point =
1000, ANN-DisPINN]{\label{fig:loss_ANN_1000}\includegraphics[width=.50\linewidth]{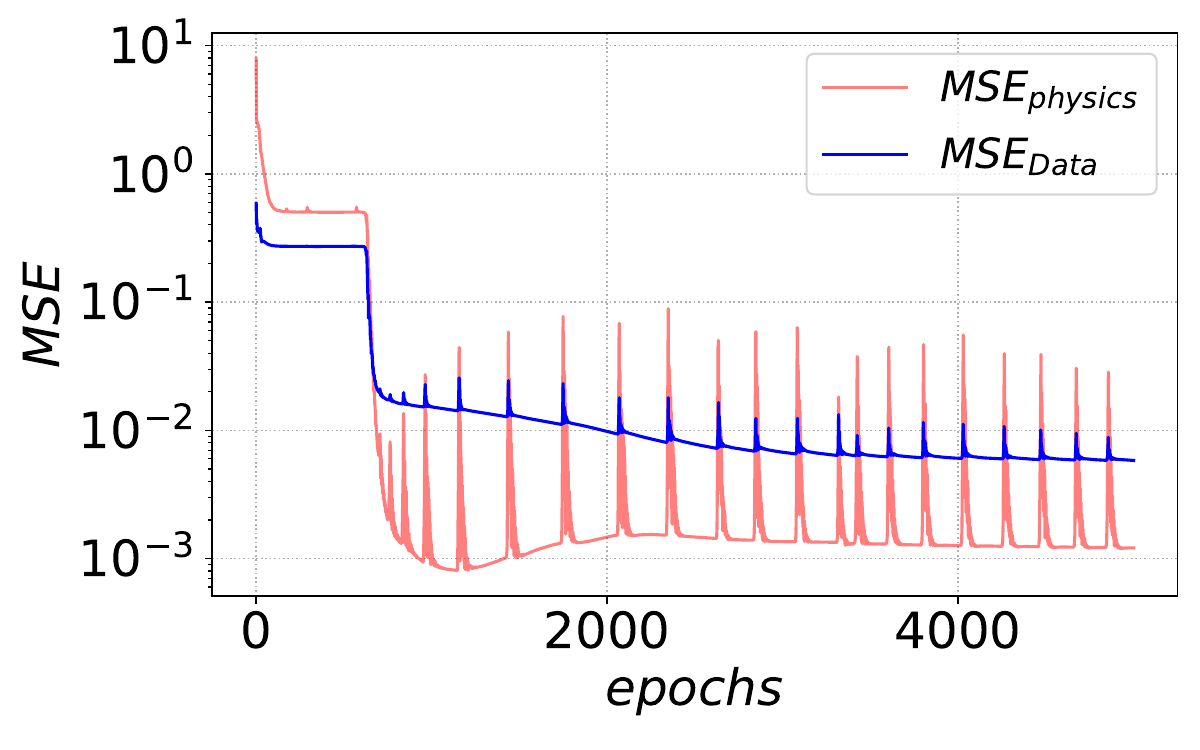}}
\end{subfigure}
\begin{subfigure}[b]{0.80\textwidth}
\centering
\subfloat[ Data Point = 1, LSTM-DisPINN]{\label{fig:loss_LSTM_1}\includegraphics[width=.50\linewidth]{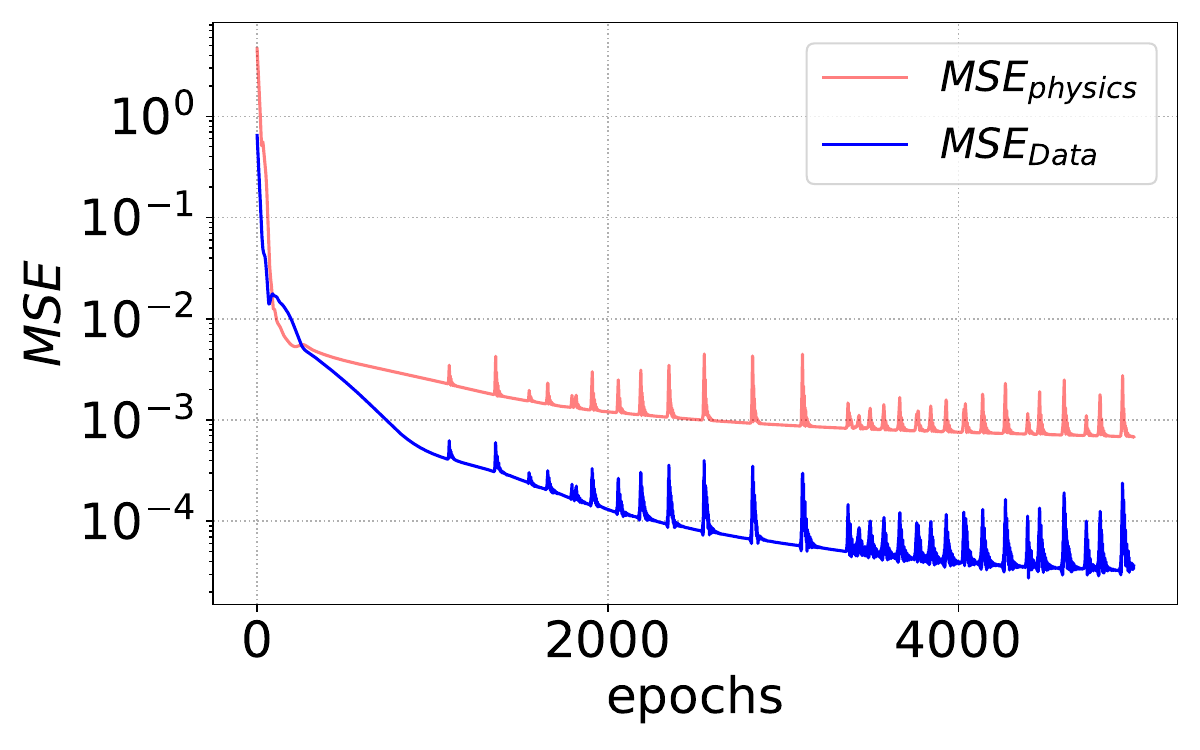}}
\subfloat[Data Point =
1000, LSTM-DisPINN]{\label{fig:loss_LSTM_1000}\includegraphics[width=.50\linewidth]{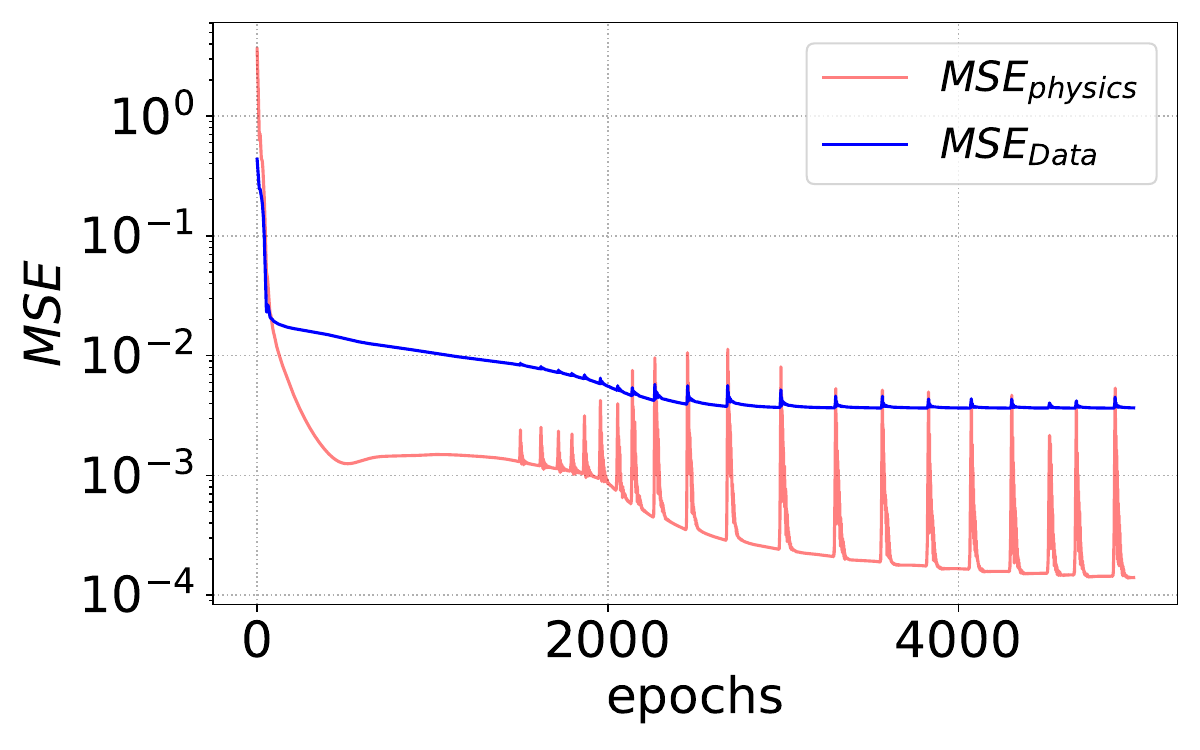}}
\end{subfigure}
\medskip
\caption{Data-driven and physics-driven loss term with no. of epochs for different physics-informed neural networks.}
\label{fig:PINN-MASSSPRING_loss}
\end{figure}

To develop the high-fidelity numerical model and surrogate model for the pitch-plunge system, the structural properties, mentioned in \autoref{eq:9} are taken the same as those defined in \cite{karnick_venkatraman_2017,halder2023deep}. 
First, the dynamical system is excited with a sinusoidal input, $\sin(10t)$ as $C_l$ coefficient and $\sin(20t)$ as $C_m$ coefficient, whereas, the coefficient $(V^{*2}/\pi)$ describing the non-dimensionalized uniform air velocity is taken as $1$. The time derivatives are computed with respect to the non-dimensional time $\tau$, where $\tau = {\omega}_\alpha t $. $t$ is the dimensional time. In the current computation, the dimensional and non-dimensional time steps are $\Delta t = \pi/1800$ and $\Delta \tau = \pi/18$ respectively. The time derivatives are discretized using the second-order backward Euler scheme. 
Next, to develop the surrogate model, an ANN and an LSTM network are generated to map between the input which
is $C_l$ and $C_m$ coefficient and the output, $h$ and $\alpha$. $1000$ time steps are considered for developing the neural networks. For the ANN network, $4$ hidden layers are considered where the layers consist of $124$, $64$, $24$ and $8$ neurons respectively. The LSTM network, on the other hand, consists of a layer with $10$ neurons in each LSTM cell. The sequence length of the LSTM cell is taken as $10$. Furthermore, we also carry out the prediction of structural response, pitch and plunge using conventional PINN (ANN, with AD-based derivative). Here, the time instants, $t$ need to be considered as input to compute the temporal derivatives from the computational graph, in physics-based loss term derived from \autoref{eq:9}. Temporal variation of the structural responses, $h$ and $\alpha$ at different time instants are considered as the output of the neural network. The hyper-parameters are the same as the ANN and LSTM networks, mentioned earlier. Therefore, We intend to carry out $5$ different surrogate models such as a data-driven approach - ANN and LSTM, Conventional AD-based PINN, AD-PINN and our proposed ANN-DisPINN and LSTM-DisPINN. \autoref{fig:PINN-MASSSPRING_loss} shows both physics-driven and data-driven $MSE$ loss terms with the number of epochs. 

We demonstrate in \autoref{fig:PINN-MASSSPRING_loss}, that when conventional AD-PINN is used, the decay rate of the loss terms is slower than the discretized physics-based PINN, ANN-DisPINN and LSTM-DisPINN. Choice of input often contributes to the convergence of the $MSE$ loss terms as shown by \cite{demo2023extended} where an extra forcing term in the input of the PINN accelerated the loss convergence. The $C_l$ and $C_m$ are the better choices as inputs in the case of the discretized-physics-based PINN as compared to the time, $t$ as an input for the conventional AD-based PINN, since the input, $C_l$ and $C_m$ has the same frequency content similar to output $h$ and $\alpha$. Furthermore, in the case of the discretized-physics-based PINN, the data-driven loss term of the PINN network is consistent with the numerical derivatives used in the physics-driven residual as indicated in \autoref{sec:PINN}. Therefore, the physics-based loss term ($MSE_{\text{physics}}$) and data-driven loss term ($MSE_{\text{Data}}$) associated with the disPINN decays faster than the AD-based PINN. It is evident from \autoref{fig:loss_AD_1000}, in the case of AD-PINN, that when the full dataset is used, $MSE_\text{data}$ does not decay below $10^{(-2)}$ due to the inconsistency of the physics-driven and data-driven loss term. However, with DisPINN while using the data points of 1000 as shown in \autoref{fig:loss_ANN_1000} and \autoref{fig:loss_LSTM_1000}, $MSE_{\text{Data}}$ has dropped below $10^{(-2)}$.

\begin{figure}[ht] \label{Fig1}
\centering
\begin{subfigure}[b]{0.80\textwidth}
\centering
\subfloat[ Data Point =
1]{\label{fig:PINN-MASSSPRING_ha}\includegraphics[width=.50\linewidth]{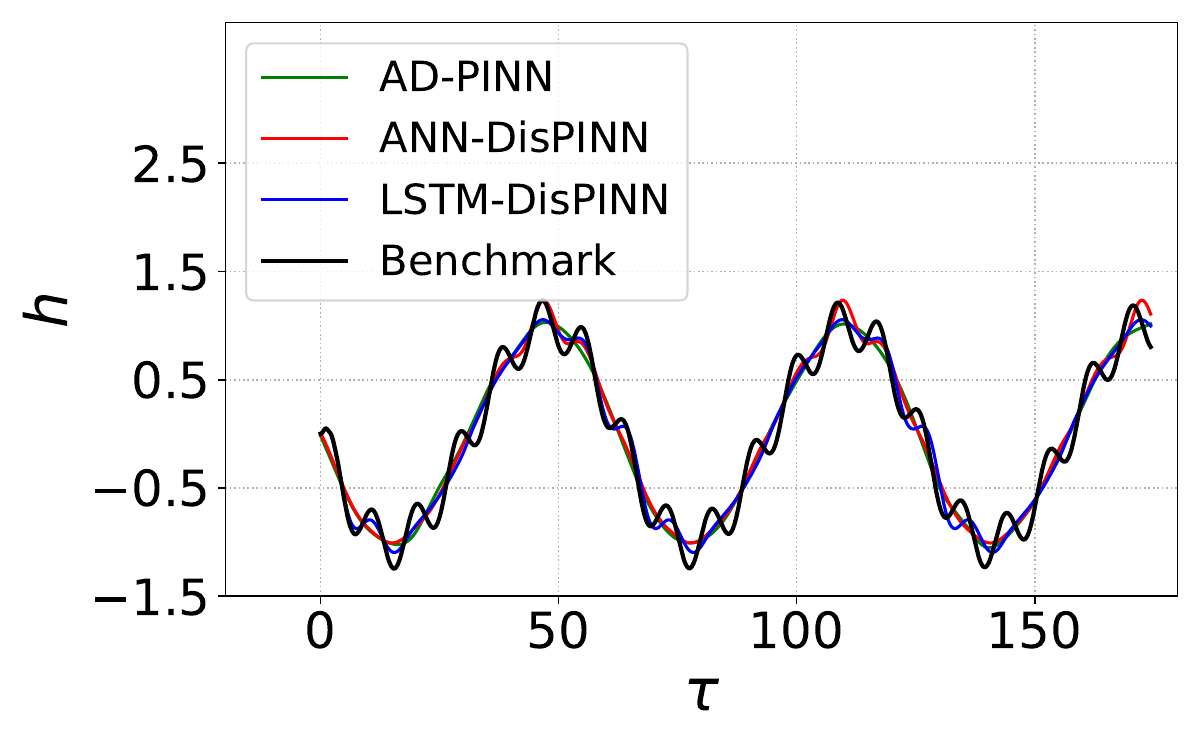}}
\subfloat[Data Point = 1
(zoomed)]{\label{fig:PINN-MASSSPRING_hb}\includegraphics[width=.50\linewidth]{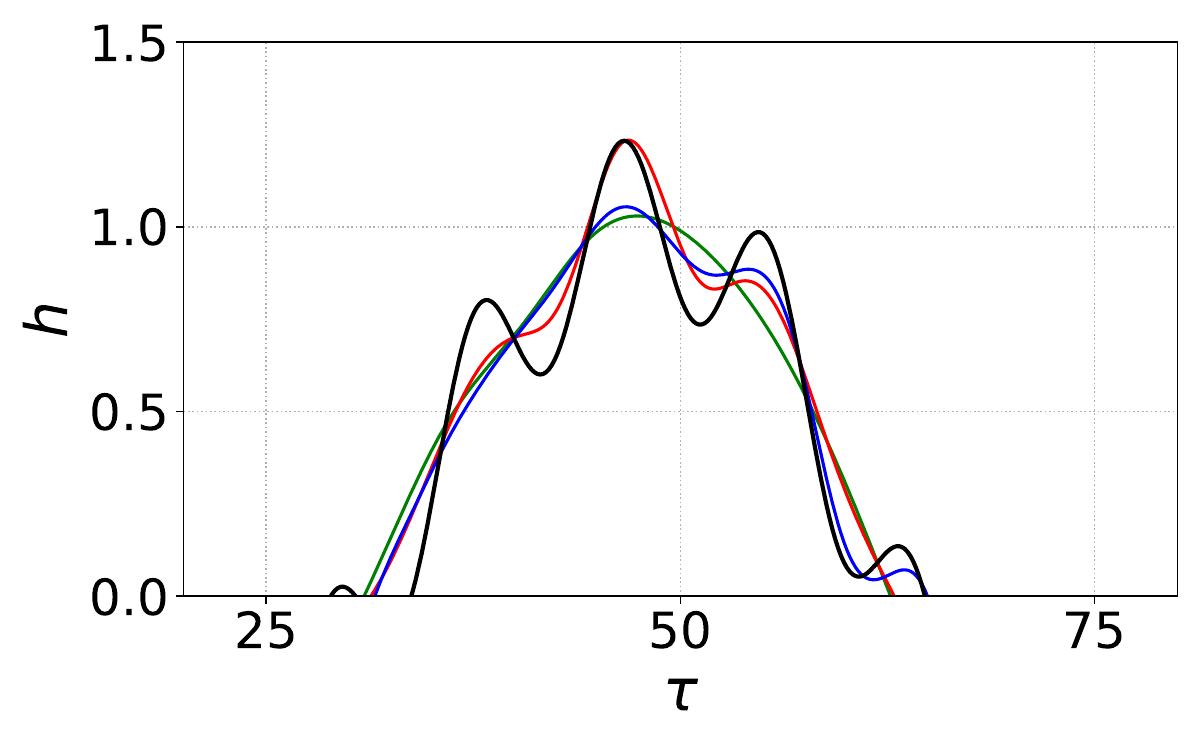}}
\end{subfigure}
\begin{subfigure}[b]{0.80\textwidth}
\centering
\subfloat[ Data Point =
1000]{\label{fig:PINN-MASSSPRING_hc}\includegraphics[width=.50\linewidth]{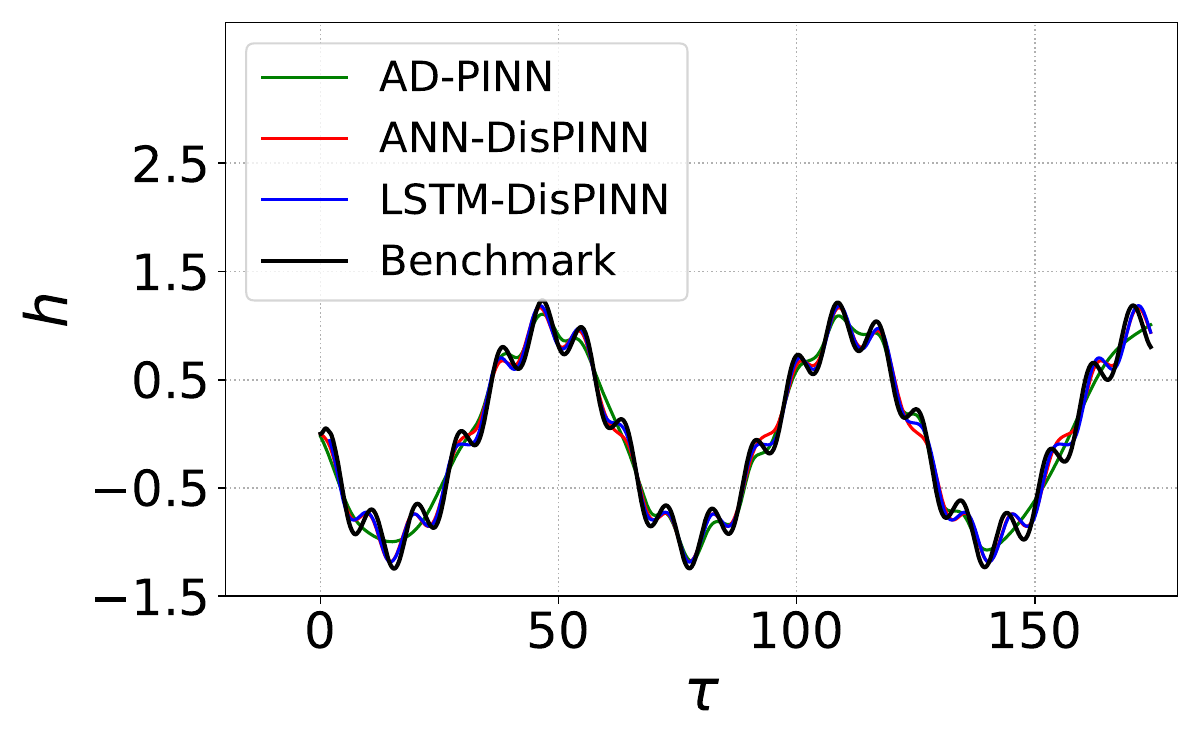}}
\subfloat[Data Point = 1000
(zoomed)]{\label{fig:PINN-MASSSPRING_hd}\includegraphics[width=.50\linewidth]{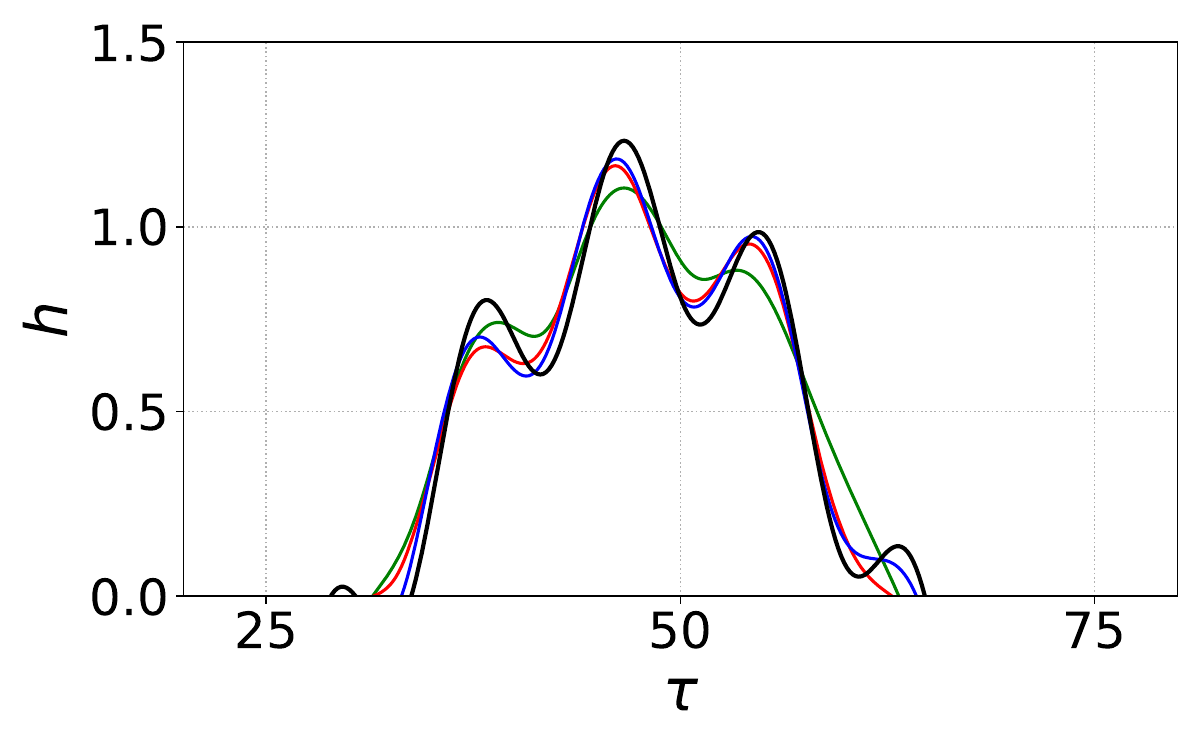}}
\end{subfigure}
\medskip
\caption{Comparison of the AD-PINN, LSTM-DisPINN and ANN-DisPINN prediction of the plunge response with benchmark numerical data}
\label{fig:PINN-MASSSPRING_h}
\end{figure}

\begin{figure}[ht]
\centering
\begin{subfigure}[b]{0.80\textwidth}
\centering
\subfloat[ Data Point =
1]{\label{fig:pitchAD1}\includegraphics[width=.50\linewidth]{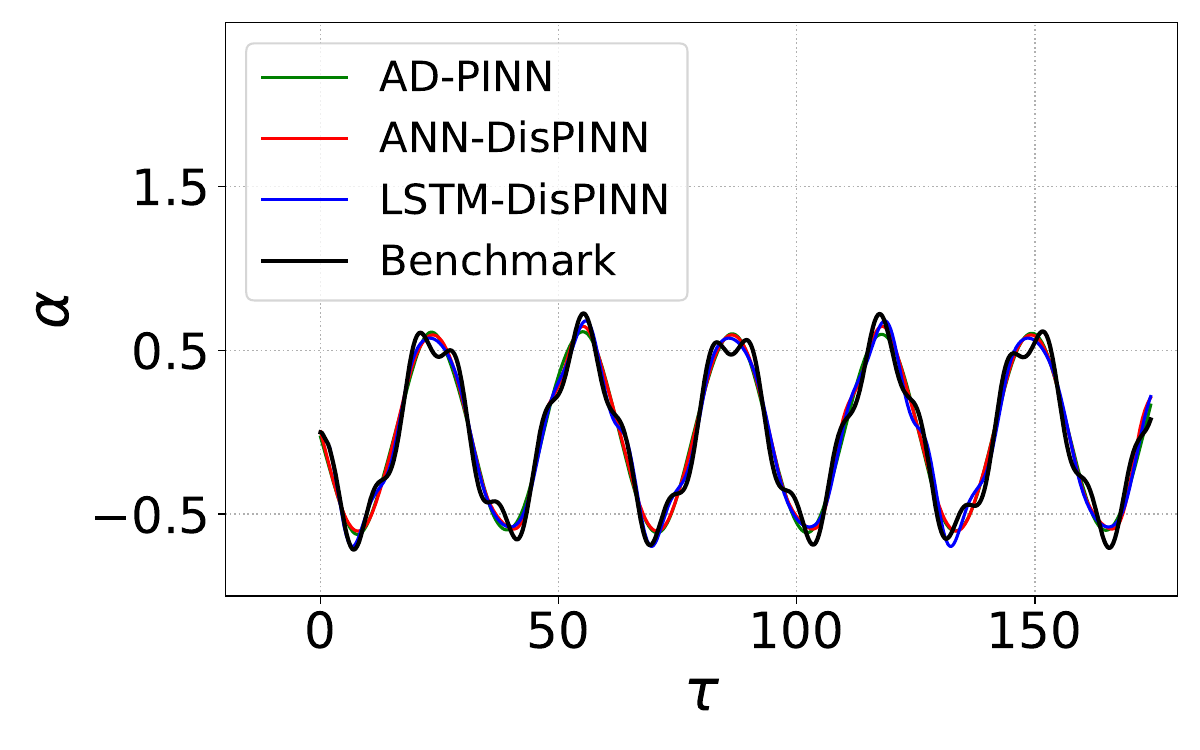}}
\subfloat[Data Point = 1
(zoomed)]{\label{fig:pitchcloseAD1}\includegraphics[width=.50\linewidth]{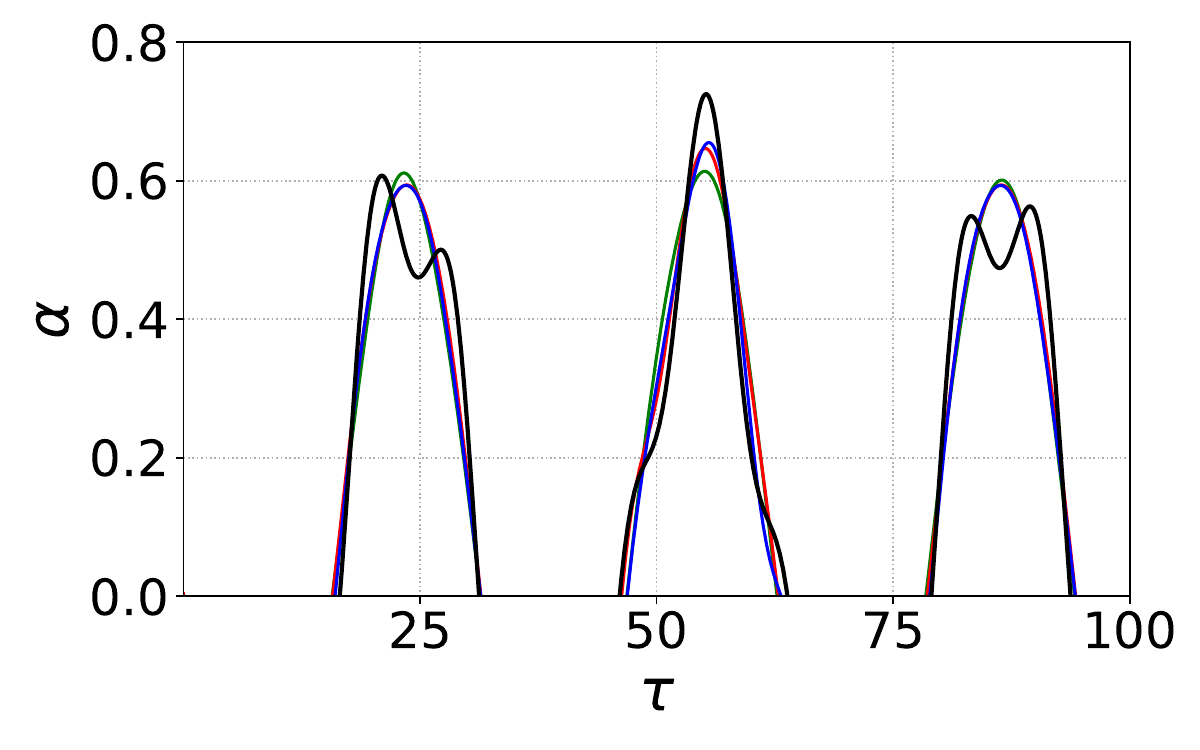}}
\end{subfigure}
\begin{subfigure}[b]{0.80\textwidth}
\centering
\subfloat[ Data Point =
1000]{\label{fig:pitchAD1000}\includegraphics[width=.50\linewidth]{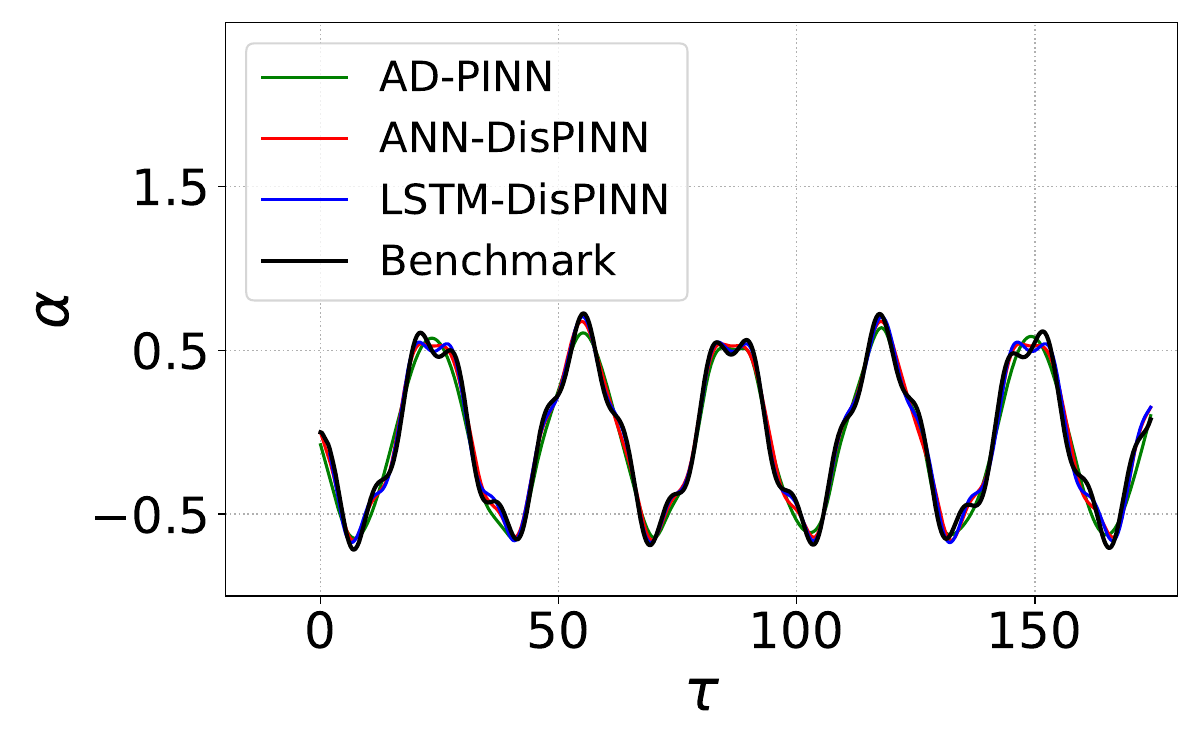}}
\subfloat[Data Point = 1000
(zoomed)]{\label{fig:pitchADclose1000}\includegraphics[width=.50\linewidth]{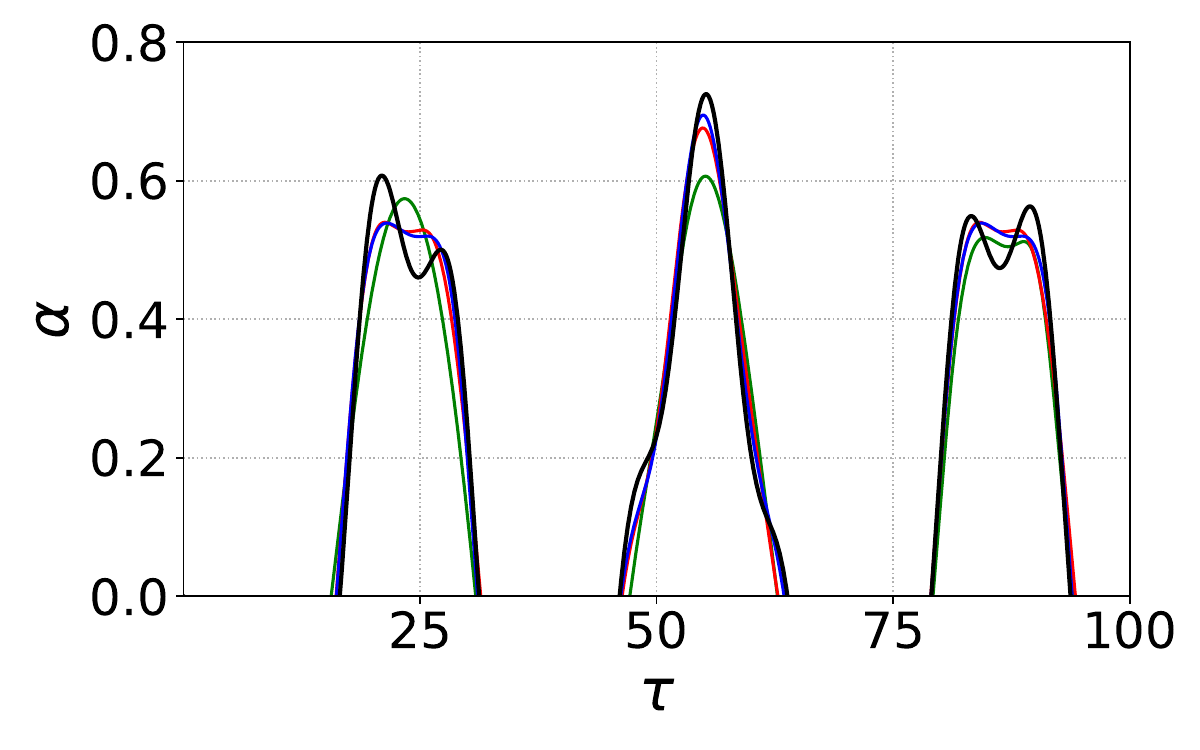}}
\end{subfigure}
\medskip
\caption{Comparison of the AD-PINN, LSTM-DisPINN and ANN-DisPINN prediction of the pitch response with benchmark numerical data}
\label{fig:PINN-MASSSPRING_a}
\end{figure}

\autoref{fig:PINN-MASSSPRING_h} and \autoref{fig:PINN-MASSSPRING_a} 
compares the prediction of the pitch and plunge response using the surrogate models i.e., the purely data-driven models such as ANN and LSTM, AD-based conventional PINN, AD-PINN and discretized physics-based PINN, ANN-DisPINN and LSTM-DisPINN networks with the high-fidelity benchmark results. The objective of the current exercise is to reconstruct the full unsteady response using only a few data points corresponding to different time instants in the structural responses. However, as mentioned in \autoref{eq:11}, the residuals from the governing equations in the discretized form (for DisPINN) and continuous form (AD-PINN) are first obtained at all the time instants in the temporal domain using the neural network prediction such as shown in \autoref{eq:19} and \autoref{eq:20} at every iteration 
of the optimization algorithm used in the neural network. Next, the residuals are used as the physics-based loss penalty term in addition to the data-driven loss term in the PINN to obtain the corrected weight matrices and bias vectors. Since there are $1000$ time instants in the temporal domain, in the ANN network, the input size consisting of sinusoidal $C_l$ and $C_m$  is $\mathbb{R}^{(1000 \times 2)}$. However, in the LSTM-DisPINN network, considering the sequence length of 10, the input tensor dimension is $\mathbb{R} ^{(990 \times 10 \times 2)}$. Furthermore, the output dimension for the ANN-DisPINN and LSTM-DisPINN is $\mathbb{R} ^{(1000 \times 2)}$ and $\mathbb{R} ^{(990 \times 2)}$ consisting of the structural responses $h$ and $\alpha$ corresponding to the input $C_l$ and $C_m$ coefficients. Here, datasets of different sizes named as data points are randomly chosen from the temporal domain. In the PINN network, out of the $1000$ training points mentioned above, only $1$, $3$, $200$, and the entire dataset, $1000$ are considered for the output prediction. These data points are different from the initial point, which are the first and $10^{th}$ time instant in the temporal domain for the ANN and LSTM network, respectively. The ANN-DisPINN and the LSTM-DisPINN are capable of predicting the system responses significantly well even when only $1$ data point is considered, as shown in \autoref{fig:PINN-MASSSPRING_ha} and zoomed view in \autoref{fig:PINN-MASSSPRING_hb}. AD-PINN, on the other hand, with fewer amount of data-points (i.e. $1$), fails to reconstruct the small oscillations in the structural response as shown in \autoref{fig:PINN-MASSSPRING_hb}. However, small oscillations in the pitch and plunge responses are well-captured with $1000$ data points using both DisPINN and AD-PINN methods as shown in \autoref{fig:PINN-MASSSPRING_hc} and zoomed view in \autoref{fig:PINN-MASSSPRING_hd}. Similar trends are also visible in the case of pitch responses as well as shown in \autoref{fig:PINN-MASSSPRING_a}.   

\begin{figure}[ht]
\centering
\begin{subfigure}[b]{0.90\textwidth}
\centering
\subfloat[  $h$ ]{\label{fig:errorms_h}\includegraphics[width=.50\linewidth]{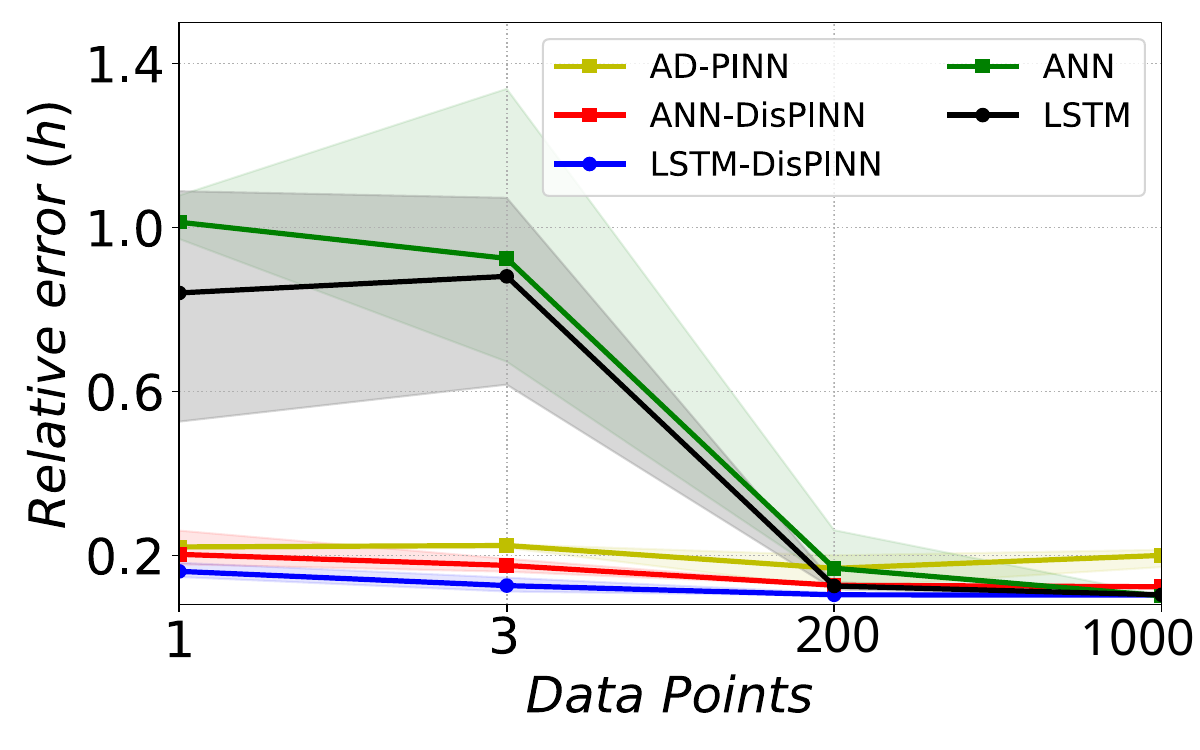}}
\subfloat[  $h$ (zoomed)]{\label{fig:errorms_h_close}\includegraphics[width=.50\linewidth]{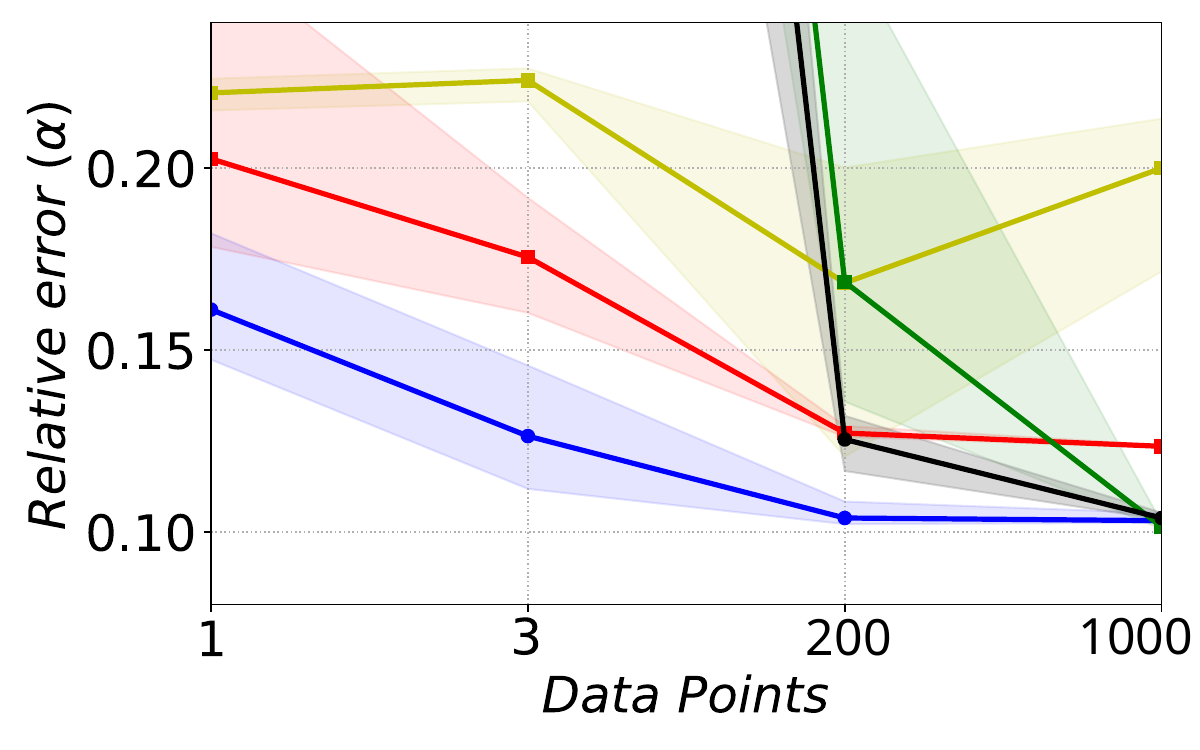}}
\end{subfigure}
\begin{subfigure}[b]{0.90\textwidth}
\centering
\subfloat[ $\alpha$ ]{\label{fig:errorms_a}\includegraphics[width=.50\linewidth]{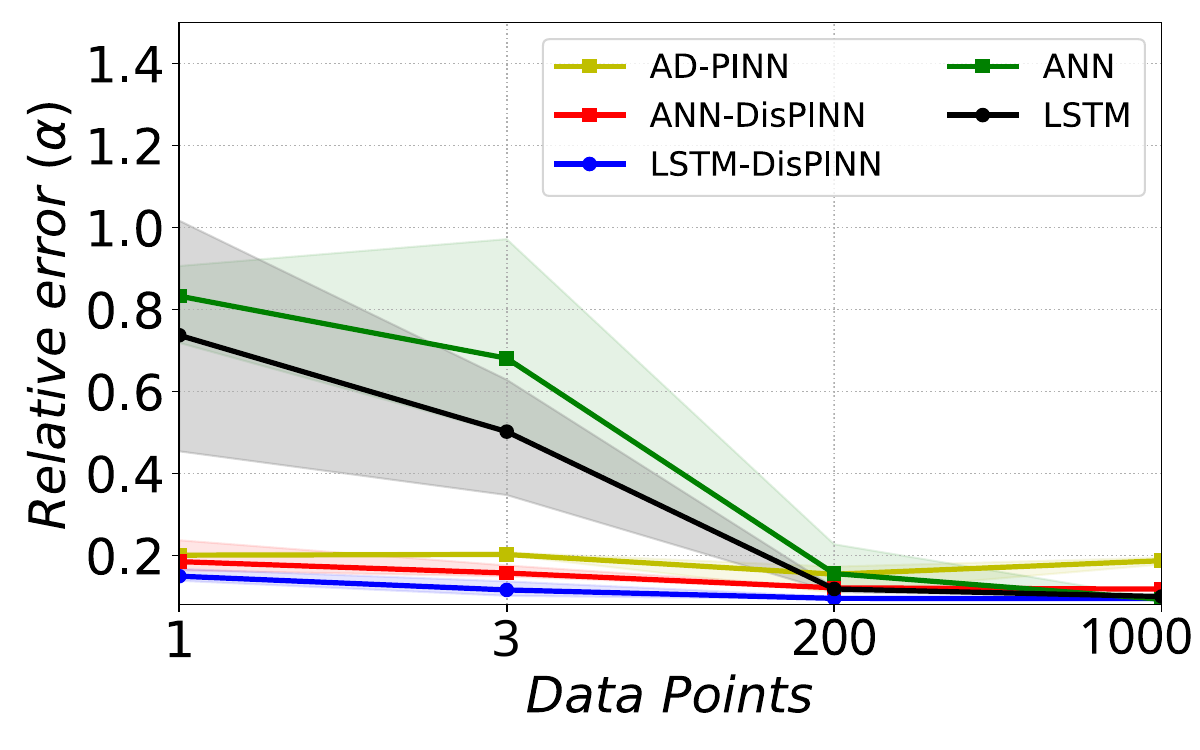}}
\subfloat[ $\alpha$ (zoomed)]{\label{fig:errorms_a_close}\includegraphics[width=.50\linewidth]{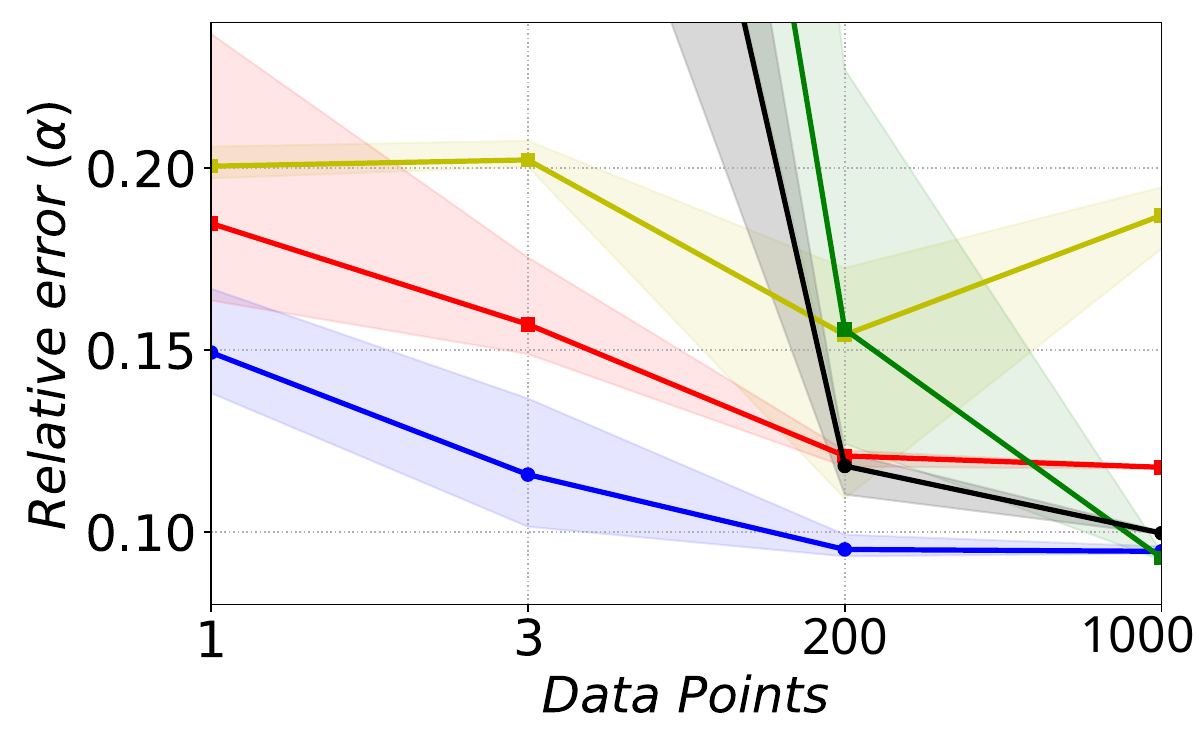}}
\end{subfigure}
\medskip
\caption{Comparison of the relative errors of prediction of the plunge and pitch response using Data-Driven and different PINN approach.}
\label{fig:error_PINN-MASSSPRING}
\end{figure}

\autoref{fig:error_PINN-MASSSPRING} shows the potential of the LSTM-DisPINN and ANN-DisPINN over pure data-driven ANN and LSTM networks and AD-PINN. The prediction errors at all the time instants using different data points are computed in the following \autoref{eq:errorpred}. 

\begin{equation}\label{eq:errorpred}
\begin{aligned}
\text{relative error} = \frac{\sqrt{\sum_{N_t}\left(y_{\text {pred}}-y_{\text {actual}}\right)^2}}{\sqrt{\sum_{N_t}\left(y_{\text {actual }}\right)^2}},
\end{aligned}
\end{equation}

where, $y_{\text {pred}}$ is the prediction from the surrogate model whereas the $y_{\text {actual}}$ is the high-fidelity results at every time instants. For the ANN network (with and without the physics constraints), $N_t = 1000$, whereas for the LSTM network (with and without the physics constraints), $N_t = 990$. We have conducted neural network training $5$ times keeping the number of data points the same at each run but changing the temporal locations of the data points to assess the maximum, minimum and the mean relative error associated with those 5 different runs. In those $5$ different training iterations, the data point(s) are selected randomly from the temporal domain and hyper-parameters are kept the same. We compute the mean, maximum, and minimum of the relative errors corresponding to the $5$ different training runs, which are termed as $E_{mean}$, $E_{max}$ and $E_{min}$, respectively.

As noticed there, for a very sparse training dataset when $1$ and $3$ randomly chosen high-fidelity data points are available in the temporal domain, both ANN and LSTM network without any physics information, produces large $E_{mean}$, whereas, AD-PINN, ANN-DisPINN and LSTM-DisPINN produce significantly accurate results. The $E_{mean}$ associated with the plunge response corresponding to ANN and LSTM networks is $1.01$ and $0.84$ for data-point $1$ and $0.92$, $0.88$ for $3$ data points. However, when the physics-based loss term is added to the data-driven loss term, with only $1$ data point, the $E_{mean}$ decreases to $0.22$ for AD-PINN and $0.18$ and $0.16$ for ANN-DisPINN and LSTM-DisPINN. The $E_{max}$ and $E_{min}$ associated with the ANN and LSTM network show that with the changes in the positions of the data points, prediction error associated with both networks can vary over a wide range as shown in \autoref{fig:errorms_h} and \autoref{fig:errorms_h_close}. However, with the addition of the physics constraints, LSTM-DisPINN performs better than the ANN-DisPINN with $1$ and $3$ data points. The average relative error associated with AD-PINN is comparatively higher than ANN-DisPINN and LSTM-DisPINN for all the data-points. For the pitch responses, the trends are similar, as shown in \autoref{fig:errorms_a} and \autoref{fig:errorms_a_close}. However, as we increase the number of data points, the prediction accuracy associated with the data-driven networks (ANN and LSTM) is close to that associated with the discretized physics-driven neural networks. Furthermore, at $200$ data points, an almost similar trend of data points $1$ and $3$ are followed, such as LSTM-DisPINN outperforms ANN-DisPINN at all $5$ training iterations and $E_{mean}$s associated with the ANN and LSTM network are above the ones corresponding to discretized physics-constrained networks. However, $E_{mean}$ associated with the purely data-driven LSTM network is higher than the AD-PINN. Conversely, when all the datasets are used in the PINN framework, the data-driven LSTM and ANN network work slightly better than the discretized-physics-driven PINN network. A similar trend is also followed in the pitch responses case. In conclusion, these observations show that when a sparse dataset is available, a discretized-physics-informed neural network shows its potential over a purely data-driven network. On the contrary, when a large dataset is available, a purely data-driven approach is enough to predict the system dynamics. Adding the physics-based loss term increases the magnitude of the residual of the network at each iteration or epoch of the optimization algorithm. The additional computational burden arising from the physics-based residual often causes the loss value to get stuck at an intermediate level if the optimization algorithm in the network is not properly tuned. Hence, the pure data-driven neural network may outperform the physics-constrained neural network when a large dataset is available. The authors like to mention that the error shown here corresponds to the hyper-parameters mentioned previously. Therefore, the performances of different networks may alter with different hyper-parameters.  

\begin{figure}[ht]
\centering
\begin{subfigure}[b]{0.90\textwidth}
\centering
\subfloat[$h$, Eqns. 5 pts]{\label{fig:pred_h1}\includegraphics[width=.33\linewidth]{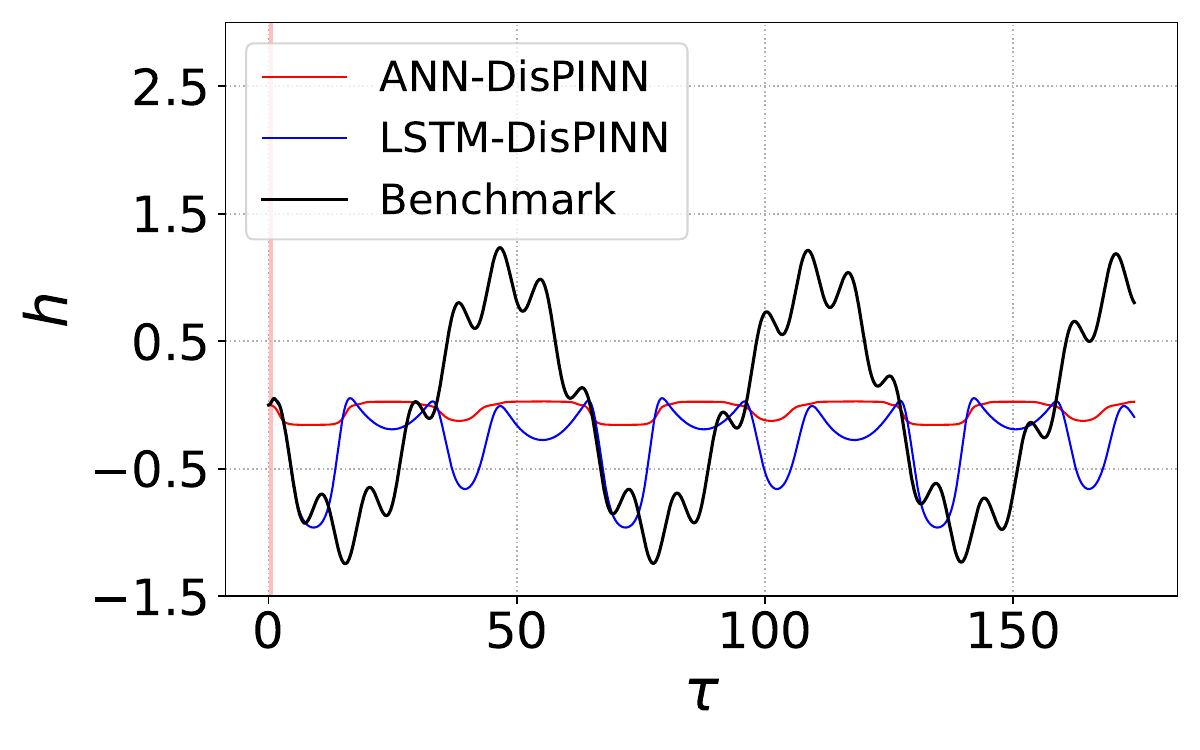}}
\subfloat[$h$, Eqns. 50 pts]{\label{fig:pred_h2}\includegraphics[width=.33\linewidth]{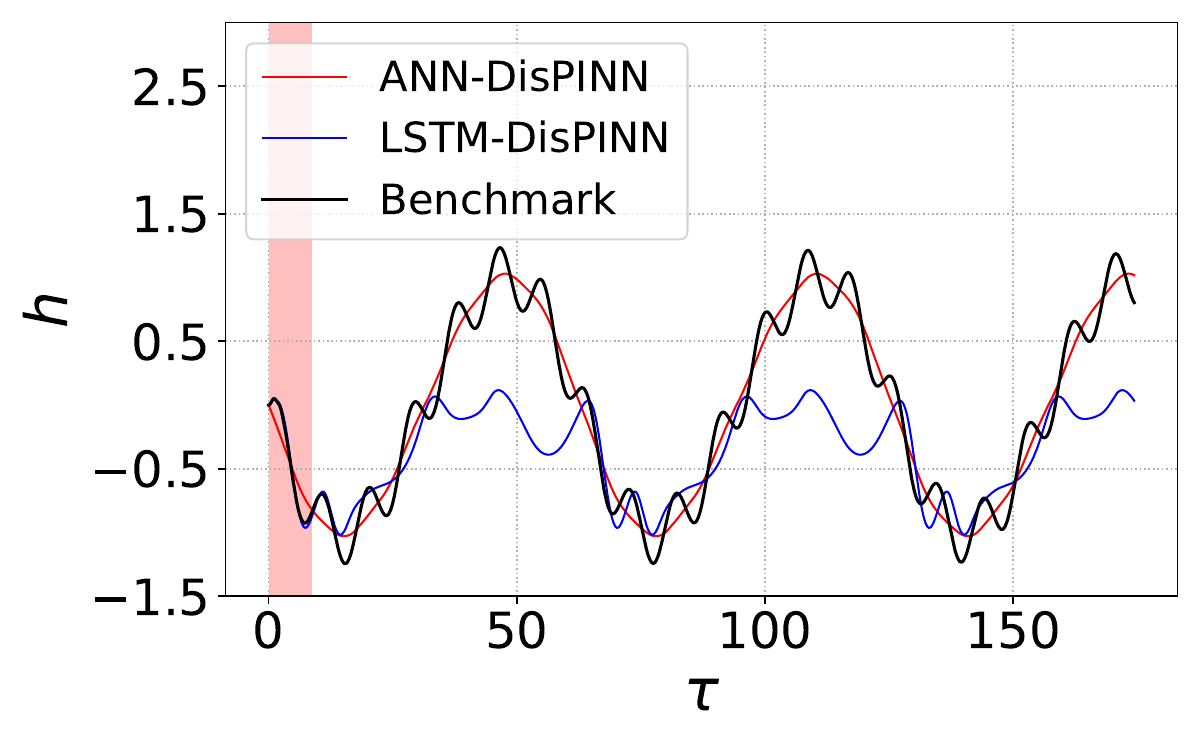}}
\subfloat[$h$, Eqns. 500 pts]{\label{fig:pred_h3}\includegraphics[width=.33\linewidth]{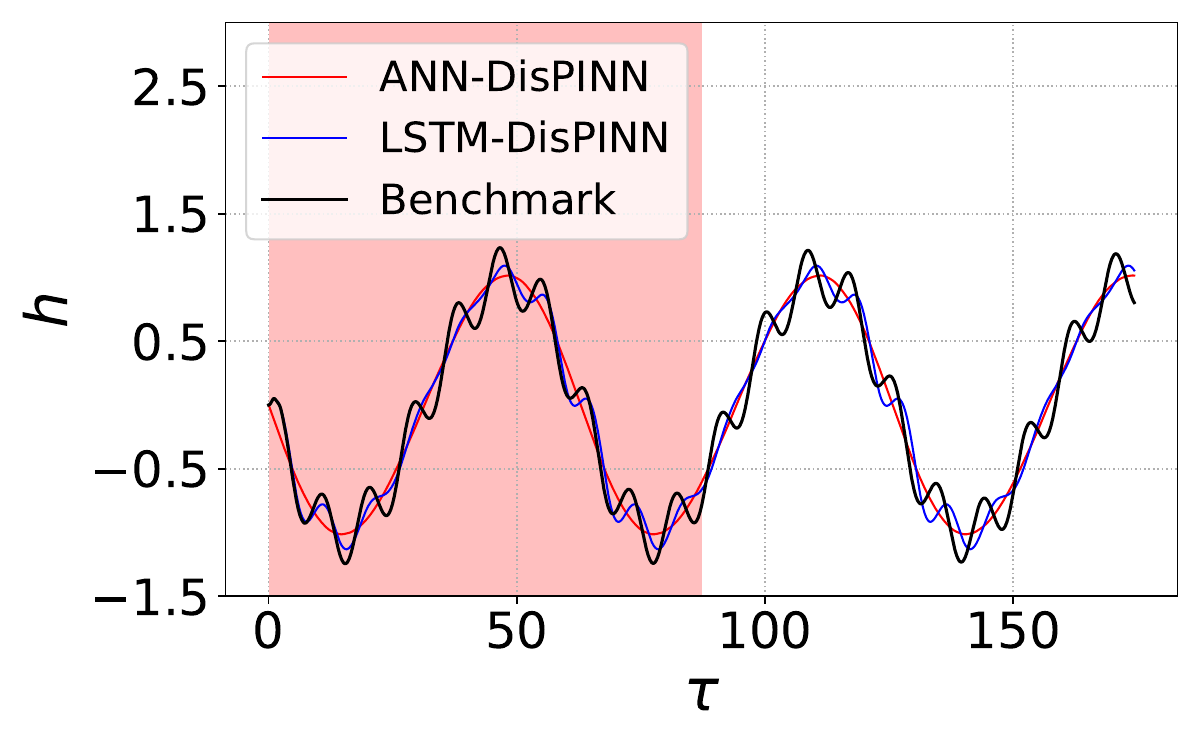}}
\end{subfigure}

\begin{subfigure}[b]{0.90\textwidth}
\centering
\subfloat[$\alpha$, Eqns. 5 pts]{\label{fig:pred_a1}\includegraphics[width=.33\linewidth]{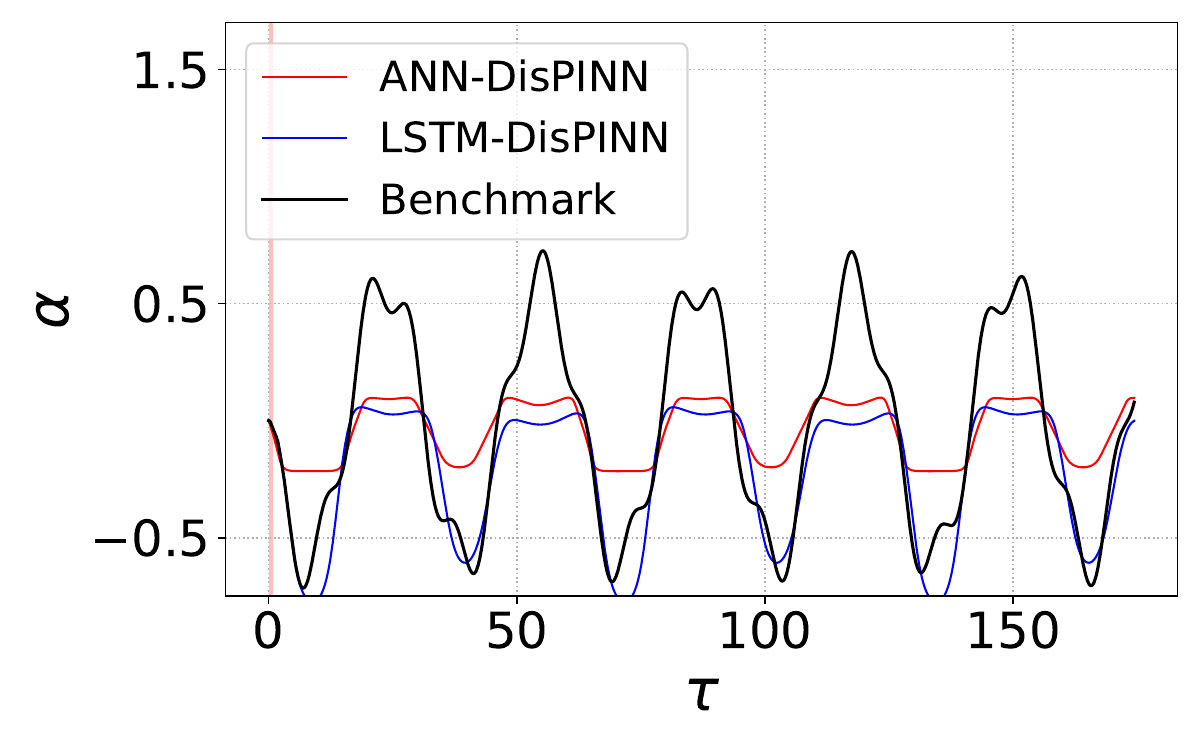}}
\subfloat[$\alpha$, Eqns. 50 pts]{\label{fig:pred_a2}\includegraphics[width=.33\linewidth]{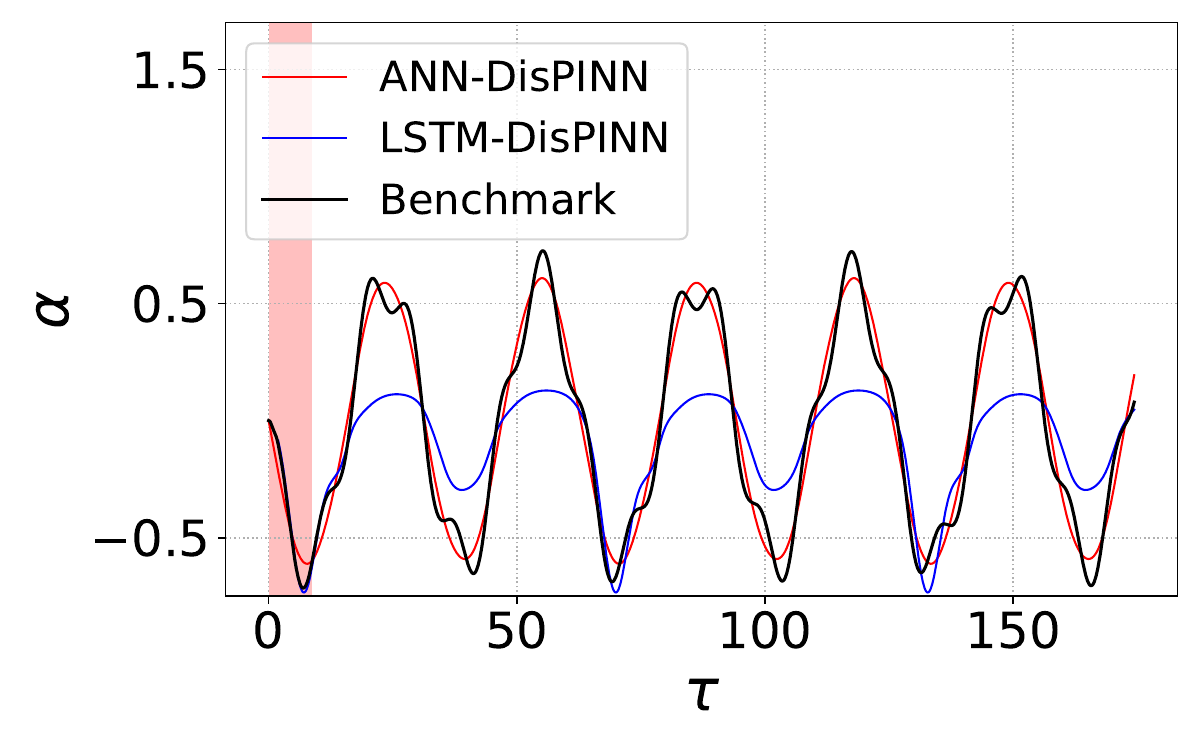}}
\subfloat[$\alpha$, Eqns. 500 pts]{\label{fig:pred_a3}\includegraphics[width=.33\linewidth]{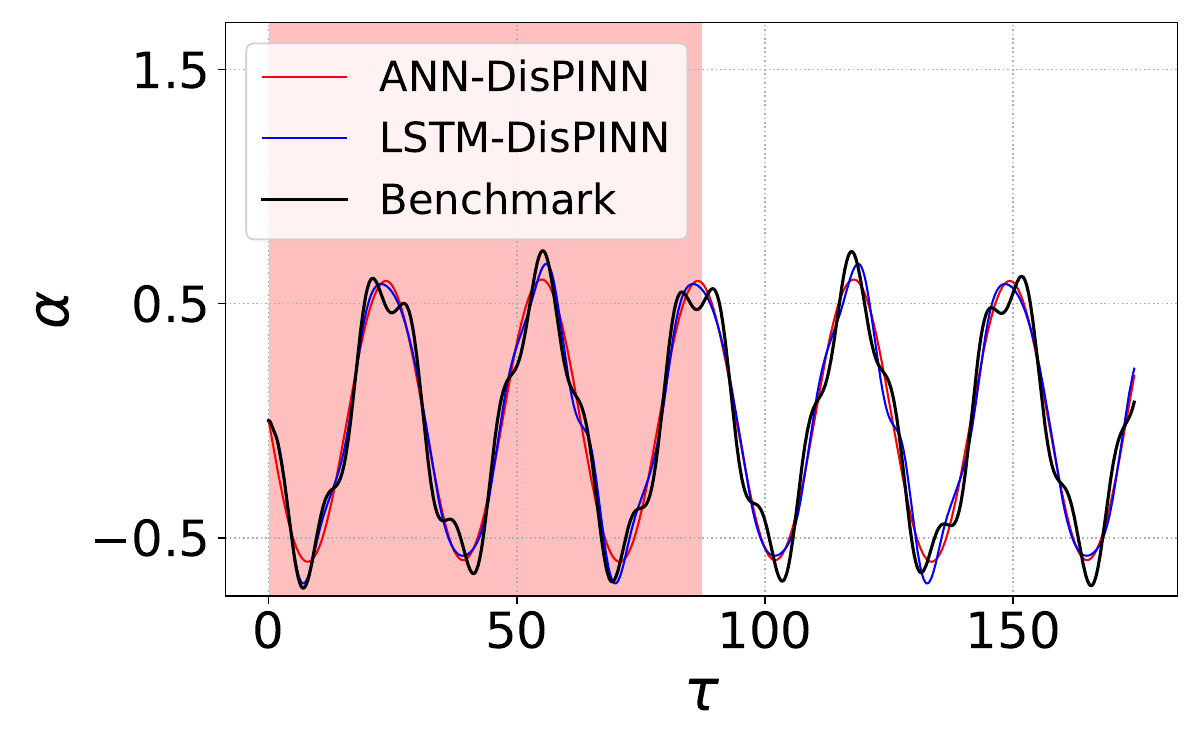}}
\end{subfigure}
\medskip
\caption{Comparison of the LSTM-DisPINN and ANN-DisPINN prediction of the plunge and pitch response with benchmark numerical data where equations solved in internal spatiotemporal points and only initial data (=1) used.}
\label{fig:pred_PINN-MASSSPRING}
\end{figure}

So far, we have demonstrated the potential of the discretized governing equation-based PINN for the reconstruction of a time series with available numerical data points at selected locations of the temporal domain, whereas the residual arising from the discretized governing equations are minimized at all the time instants. To check such networks' prediction capability, governing equations are solved at initial time instants in addition to numerical data points selected within the initial time instants. The pitch and plunge responses are now predicted at the future time instants. \autoref{fig:pred_PINN-MASSSPRING} shows the prediction capability of the ANN-DisPINN and LSTM-DisPINN network and compares them with the benchmark high-fidelity data. \autoref{fig:pred_h1}, \autoref{fig:pred_h2} and \autoref{fig:pred_h3} show the prediction capability when first $5$, $50$, and $500$ time-instants are used for the residual-minimization arising from the discretized equation and only the first time instant in the temporal domain is considered for the data-driven loss term. No additional data points are selected other than the data-driven constraints arising from the initial condition.        

As shown in \autoref{fig:pred_h1}, with $5$ points used for discretized equation, both the ANN-DisPINN and LSTM-DisPINN show the large relative errors of $0.97$ and $1.07$ as indicated in \autoref{Table 1}. However, as we increase the training points where the equations are solved in 50 temporal instants, shown in \autoref{fig:pred_h2}, the ANN-DisPINN can predict the trend of the pitch and plunge responses generating an error of $0.20$ whereas the error associated with the LSTM-DisPINN is still very large, $0.78$. Conversely, at $500$ points, as noticed in \autoref{fig:pred_h3}, the error associated with the ANN network does not change significantly, whereas the prediction error associated with the LSTM network drops to $0.20$ demonstrating the prediction capability of discretized equation-based networks with only one numerical data point considered at the initial time instant. As shown in \autoref{fig:pred_h3}, the LSTM-DisPINN network has been able to predict the high-frequency responses, whereas the ANN-DisPINN fails.      

\begin{table}[h!] 
\centering
 \begin{tabular}{||c | c | c||} 
 \hline
 Eqns Solved & ANN-DisPINN & LSTM-DisPINN \\ [0.5ex] 
 \hline\hline
 5 & 0.97 & 1.07  \\ 
 50 & 0.20 & 0.78 \\
 500 & 0.22 & 0.20 \\ [1ex] 
 \hline
 \end{tabular}
 \caption{\label{Table 1} Prediction error associated with the ANN-DisPINN and LSTM-DisPINN for plunge response}
\end{table}

\subsection{Burgers' Equation} \label{subsec:Burgers_Full}

\begin{figure}[ht]
\centering
\begin{subfigure}[b]{0.90\textwidth}
\centering
\subfloat[ Data Point =
1]{\label{fig:AD_2}\includegraphics[width=1.0\linewidth]{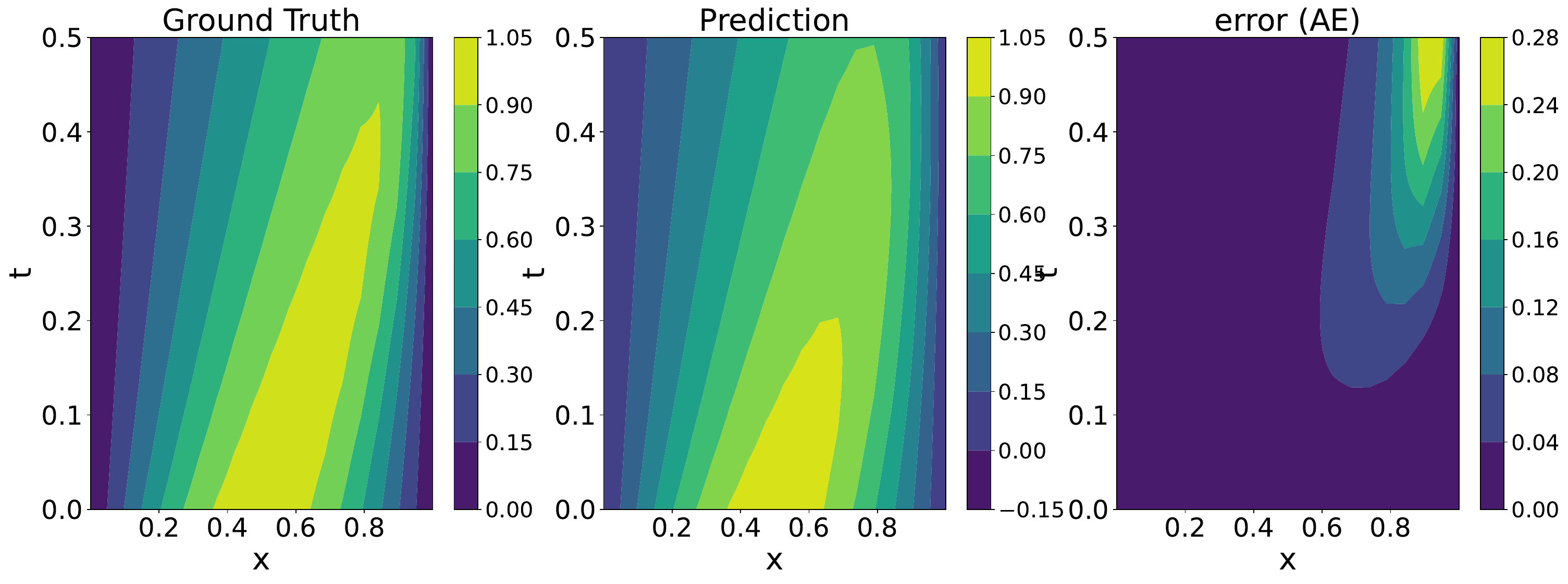}}
\end{subfigure}
\begin{subfigure}[b]{0.90\textwidth}
\subfloat[ Data Point =
1000]{\label{fig:AD_1000}\includegraphics[width=1.0\linewidth]{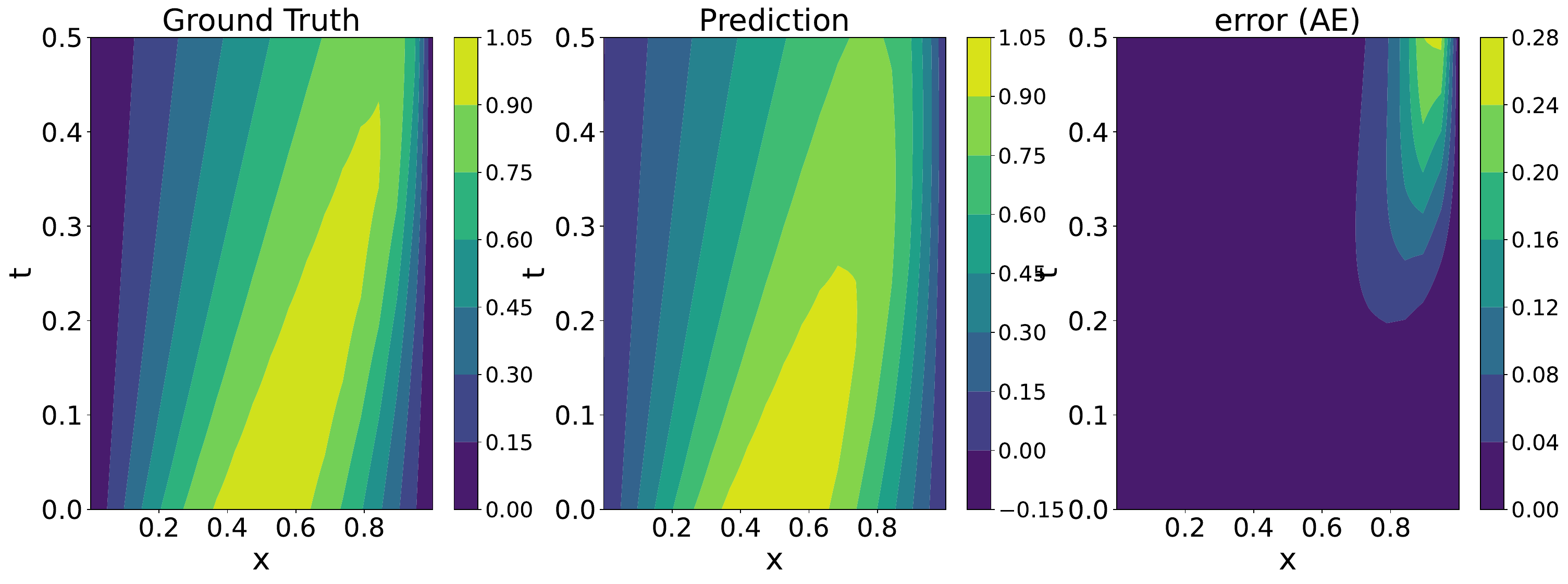}}
\end{subfigure}
\medskip
\caption{Comparison of the AD-based conventional PINN prediction with different number of data points used.}
\label{fig:AD-BurgersFull}
\end{figure}

The Burgers' equation is used as a test case to demonstrate the potential of our approaches to spatiotemporal problems. In the current subsection, the full-order discretized equation is considered in the loss term of the ANN and LSTM network in addition to the data-driven loss term. To develop the surrogate model for the Burgers' equation, $4$ layers with $124$, $64$, $24$ and $8$ neurons are considered. The ANN network's learning rate is $0.006$, and the number of epochs is $6000$. In the case of the LSTM network, similar to the mass-spring system, sequence length is considered as $10$. The number of neurons in each layer is taken as $10$. The learning rate of the LSTM network is considered $0.1$, and the number of epochs is $6000$. Furthermore, in the following \autoref{subsec:Burgers_reduced}, we will show the application of the reduced-order discretized equation obtained after the POD-Galerkin projection of the full-order system as shown in \autoref{eq:8} followed by integration with data-driven ANN and LSTM network as shown in \autoref{fig:Fig2} and \autoref{fig:Fig3}. We first assess the accuracy of the conventional PINN framework, where an ANN-based neural network is coupled with the governing equations, where the derivatives are computed using automatic differentiation. The predicted results are compared with the high-fidelity results, where the convective and diffusive parts of the equation are computed using numerical schemes demonstrated in \autoref{eq:12}. For the generation of the high-fidelity database, the temporal domain is discretized into $100$ time steps, whereas the spatial domain is discretized into a $20$ equally spaced domain. \autoref{fig:AD-BurgersFull} shows the high-fidelity results, distribution of the spatiotemporal variable arising from the surrogate model, and associated absolute error distribution between the prediction and ground truth. The automatic differentiation-based derivative requires both the space and time variables as input of the neural network to compute the gradient terms. Therefore, the size of the input of the conventional PINN is $\mathbb{R}^{2000 \times 2}$ where the size of the output variable $u$ is $\mathbb{R}^{2000 \times 1}$ which corresponds to the spatio-temporal domain ($x$ and $t$) mentioned in the input. In figure \autoref{fig:AD-BurgersFull}, different sizes of training datasets are considered i.e. the number of the data point is $1$ as shown in \autoref{fig:AD-BurgersFull}(a) and all the data points which are 1000 shown in \autoref{fig:AD-BurgersFull}(b).

\begin{figure}[ht]
\centering
\begin{subfigure}[b]{0.80\textwidth}
\centering
\subfloat[ Data Point = 1
]{\label{fig:DataANN1}\includegraphics[width=1.0\linewidth]{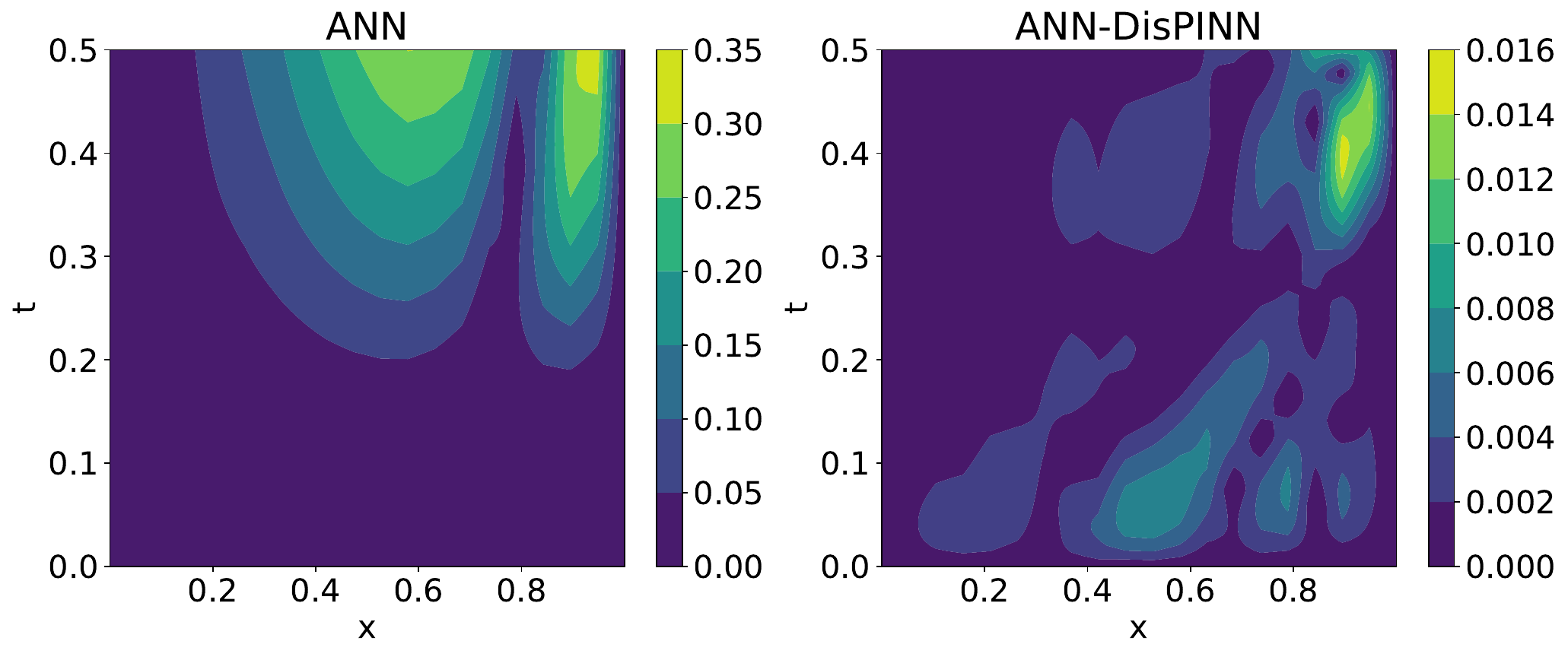}}
\end{subfigure}
\begin{subfigure}[b]{0.80\textwidth}
\centering
\subfloat[ Data Point =
1000]{\label{fig:DataANN1000}\includegraphics[width=1.0\linewidth]{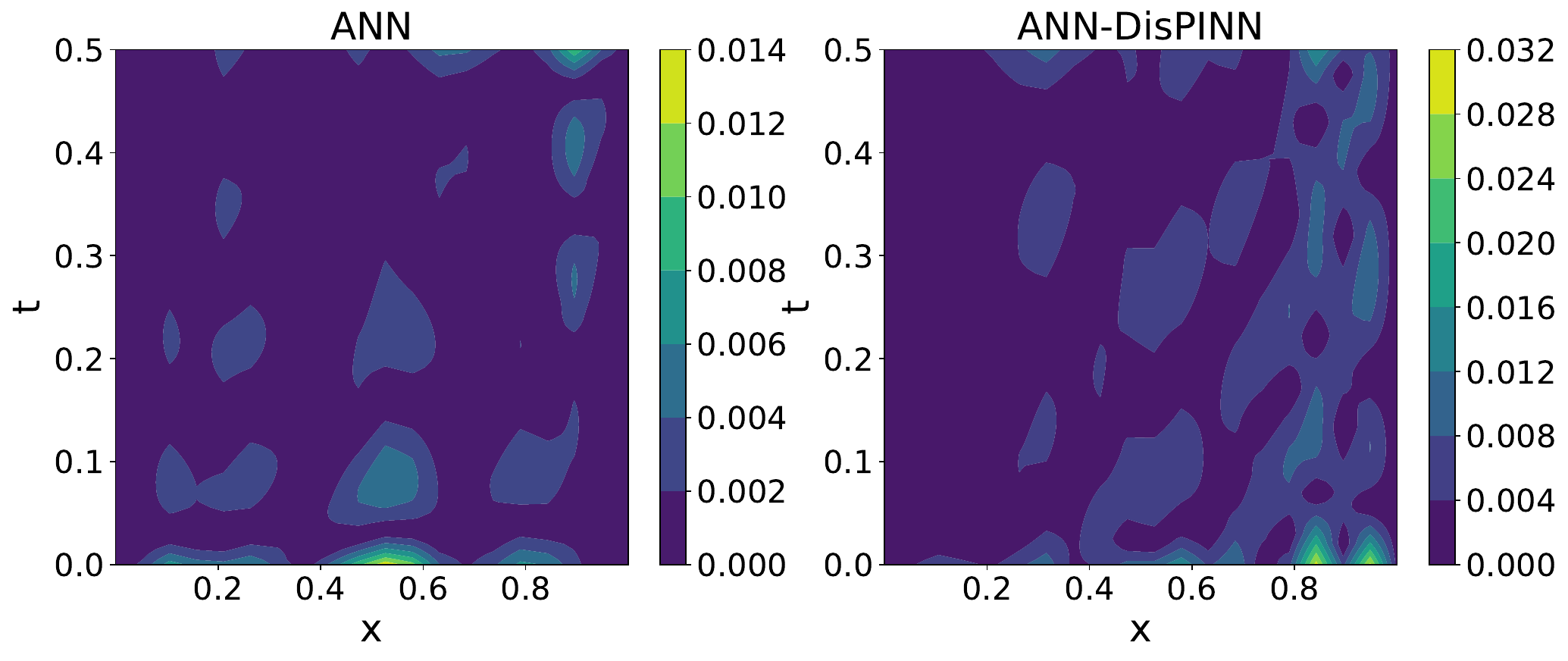}}
\end{subfigure}
\medskip
\caption{Absolute error associated with the ANN and ANN-DisPINN with different Data points used.}
\label{fig:ANN-error}
\end{figure}

\begin{figure}[ht]
\centering
\begin{subfigure}[b]{0.80\textwidth}
\centering
\subfloat[ Data Point = 1
]{\label{fig:LSTM-error1}\includegraphics[width=1.0\linewidth]{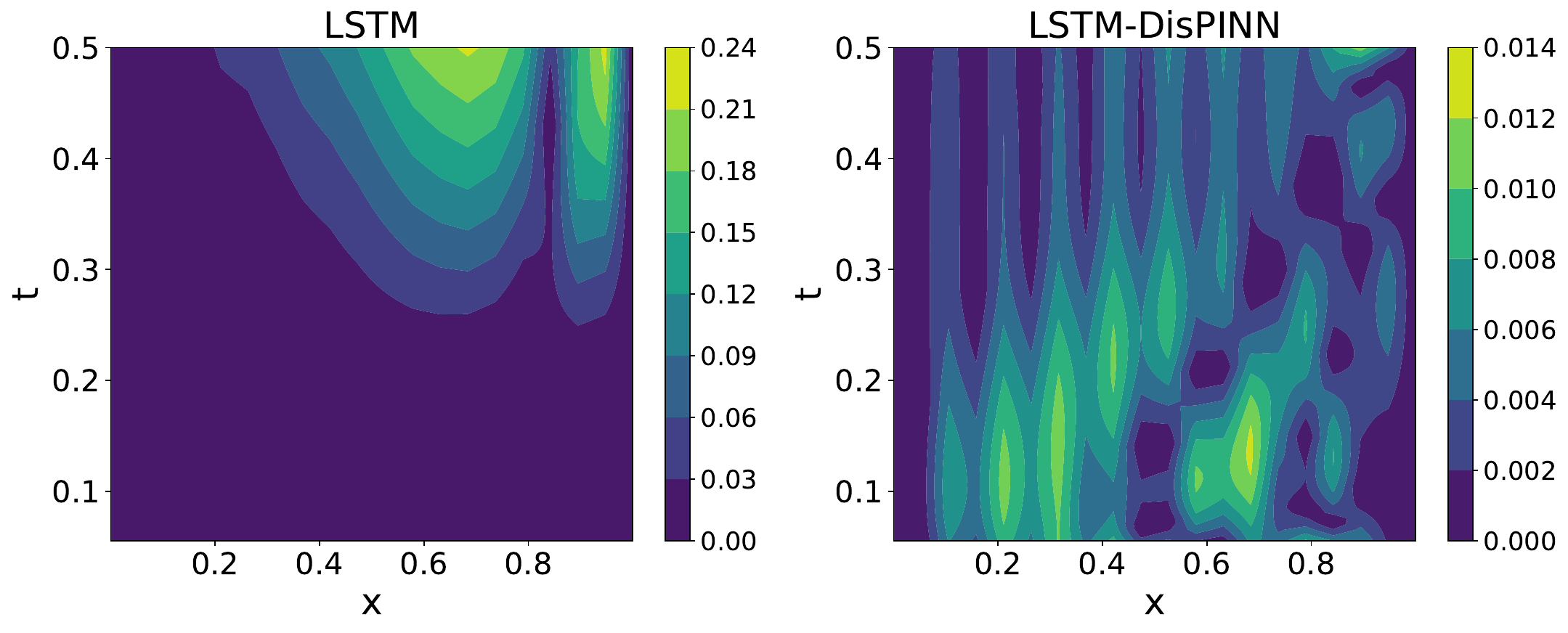}}
\end{subfigure}
\begin{subfigure}[b]{0.80\textwidth}
\centering
\subfloat[ Data Point =
1000]{\label{fig:LSTM-error2}\includegraphics[width=1.0\linewidth]{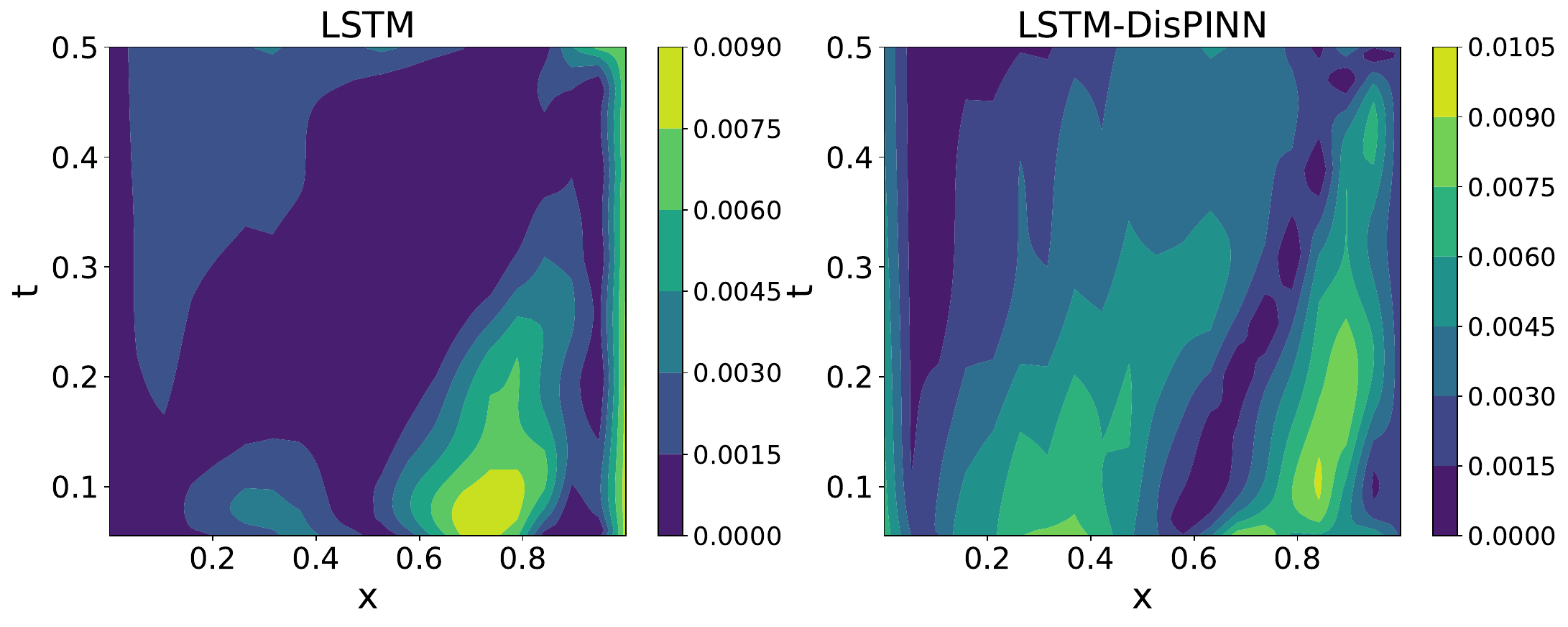}}
\end{subfigure}
\medskip
\caption{Absolute error associated with the LSTM and LSTM-DisPINN with different Data points used.}
\label{fig:LSTM-error}
\end{figure}

The data points in \autoref{fig:AD-BurgersFull}(a) are chosen randomly from the spatio-temporal domain. Apart from the data-driven loss term arising from the data point, the initial conditions, $u(x,0) = 0$, and boundary conditions, $u(0,t) = u(L,t) = 0$ are enforced in the data-driven loss part of the PINN network. However, the physics-based loss terms arising from governing equations are computed at 100 random locations in the spatio-temporal domain. The benchmark or the high-fidelity results are computed using a numerical scheme where the computation of the derivative is very different from the one adopted in the PINN network (AD-based derivative). Therefore, the data-driven loss term does not comply with the physics-based loss term thereby resulting in a large deviation of the surrogate model with respect to the high-fidelity model as shown at the end of the $x$-location and $t$-domain or the north-right corner of the error (AE) of \autoref{fig:AD-BurgersFull}(a) and \autoref{fig:AD-BurgersFull}(b). The maximum error associated with the conventional PINN, when only $1$ and full datasets are available, are respectively $0.2726$ and $0.25$ as indicated in \autoref{Table 2}. However, gain in accuracy is not achieved much as we increase the number of data points in the training set because the major source of
 error is the no-compatibility of the AD-based derivative and numerical derivative, which produced the dataset initially and not the sparse data point. Now, we introduce the application of discretized governing equation-based loss function in the PINN network with the ANN and LSTM network. As mentioned previously, conventional PINN computes the derivative using automatic differentiation, which requires the spatio-temporal coordinates $(x \ \ \text{and} \ \ t)$ as inputs. However, our approach of using discretized governing equation-based PINN offers the flexibility of choice of network input. For example, in the current work, only the time instants, i.e., $100$ time steps are considered as the network input. The $u$ variables at 20 grid points are considered as the output. Therefore for the ANN network, the input size is $\mathbb{R}^{(100 \times 1)}$ whereas the output dimension is $\mathbb{R}^{(100 \times 20)}$. In \autoref{fig:ANN-error} and \autoref{fig:LSTM-error} the error
associated with ANN-DisPINN and LSTM-DisPINN is compared with a comprehensive
data-driven ANN and LSTM network. 

The error associated with discretized governing equation-based ANN-DisPINN is much lower than the conventional
PINN network as shown in \autoref{fig:AD-BurgersFull}. Furthermore, the training time is greatly reduced since the temporal domain is only considered as an input instead of both space ($x$) and time ($t$) variables. It shows, similar to the Mass-spring system, the potential of the ANN-DisPINN and LSTM-DisPINN as compared to
the complete data-driven ANN and LSTM ROM when a very sparse dataset is
available. When one training data point is only available, the maximum error
associated with the data-driven ANN network is $0.35$ whereas, with the ANN-DisPINN
network with the discretized equation is $0.016$, indicated in \autoref{Table 2}. However, when all the data
points are used, the maximum error incurred with the ANN network is $0.014$,
whereas with the ANN-DisPINN network is $0.032$ and 
 thereby
supporting the conclusion we made in the earlier section on the mass-spring
system - when a large dataset is available, purely data-driven ROM tends to
outperform the disPINN or AD-PINN. While using the LSTM network, a sequence
length of $10$ and neuron size of $10$ in each LSTM cell is used. As shown in the
\autoref{fig:Fig3}, input is first reshaped into the sequence windows to get the tensor
structure, and the length of the spatio-temporal series becomes $90$, corresponding to the $90$ output. In the case of one data point, as shown in \autoref{fig:LSTM-error1}, the maximum error associated
with the data-driven LSTM and discretized-physics-driven LSTM networks is $0.2154$ and $0.011$. However, when all the
data points are used, both the error is $\approx0.008$. The errors are indicated in \autoref{Table 2}.

\begin{table}[h!] 
\centering
 \begin{tabular}{||c | c | c | c |c |c ||} 
 \hline
 Data-Points & AD-PINN & ANN-DisPINN & LSTM-DisPINN & ANN & LSTM \\ [0.5ex] 
 \hline\hline
 1 & 0.27 & 0.016 & 0.011 & 0.35 & 0.2154 \\ 
 1000 & 0.25 & 0.032 & 0.008 & 0.014 & 0.008\\ [1ex] 
 \hline
 \end{tabular}
 \caption{\label{Table 2} Maximum absolute error associated with the ANN-DisPINN and LSTM-DisPINN for Burgers' equation.}
\end{table}

\begin{figure}[ht]
\centering
\begin{subfigure}[b]{0.80\textwidth}
\centering
\subfloat[ ANN-DisPINN
]{\label{fig:pred_burgers1}\includegraphics[width=1.0\linewidth]{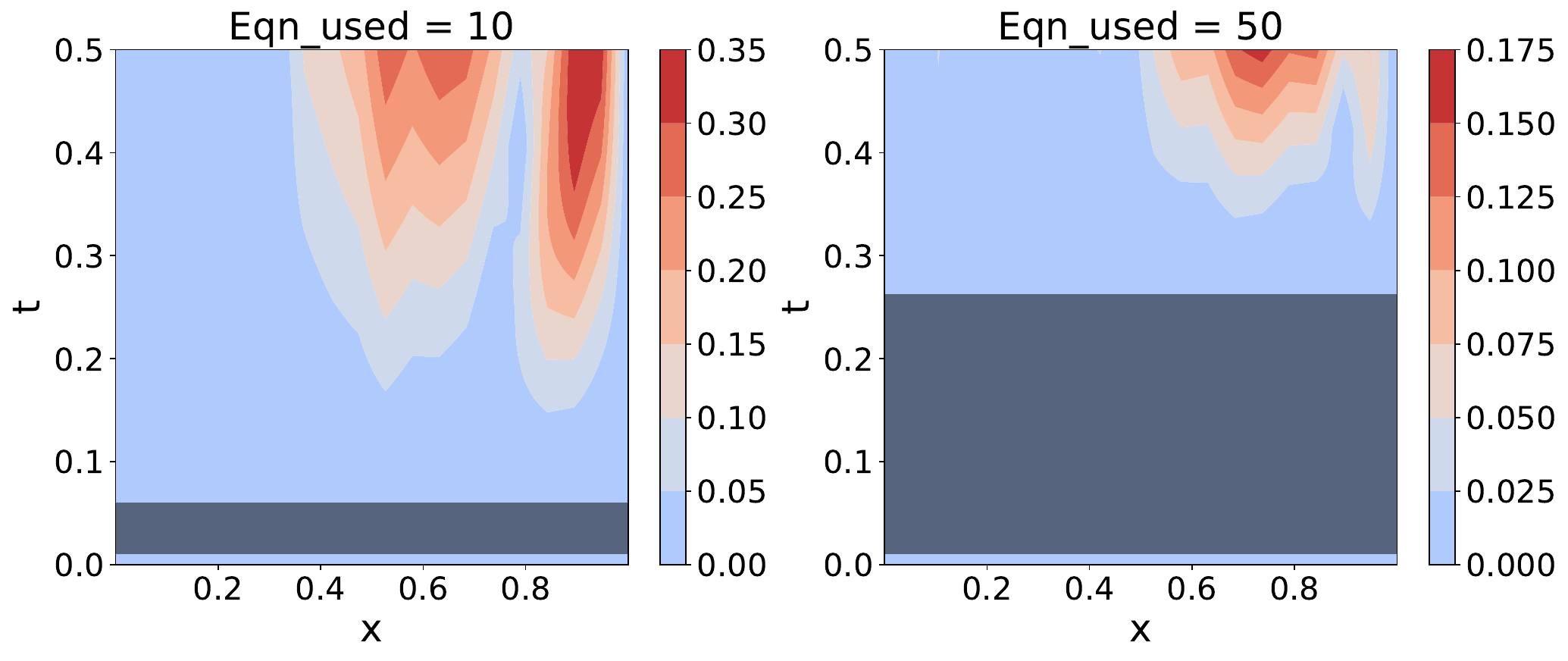}}
\end{subfigure}
\begin{subfigure}[b]{0.80\textwidth}
\centering
\subfloat[ LSTM-DisPINN]{\label{fig:pred_burgers2}\includegraphics[width=1.0\linewidth]{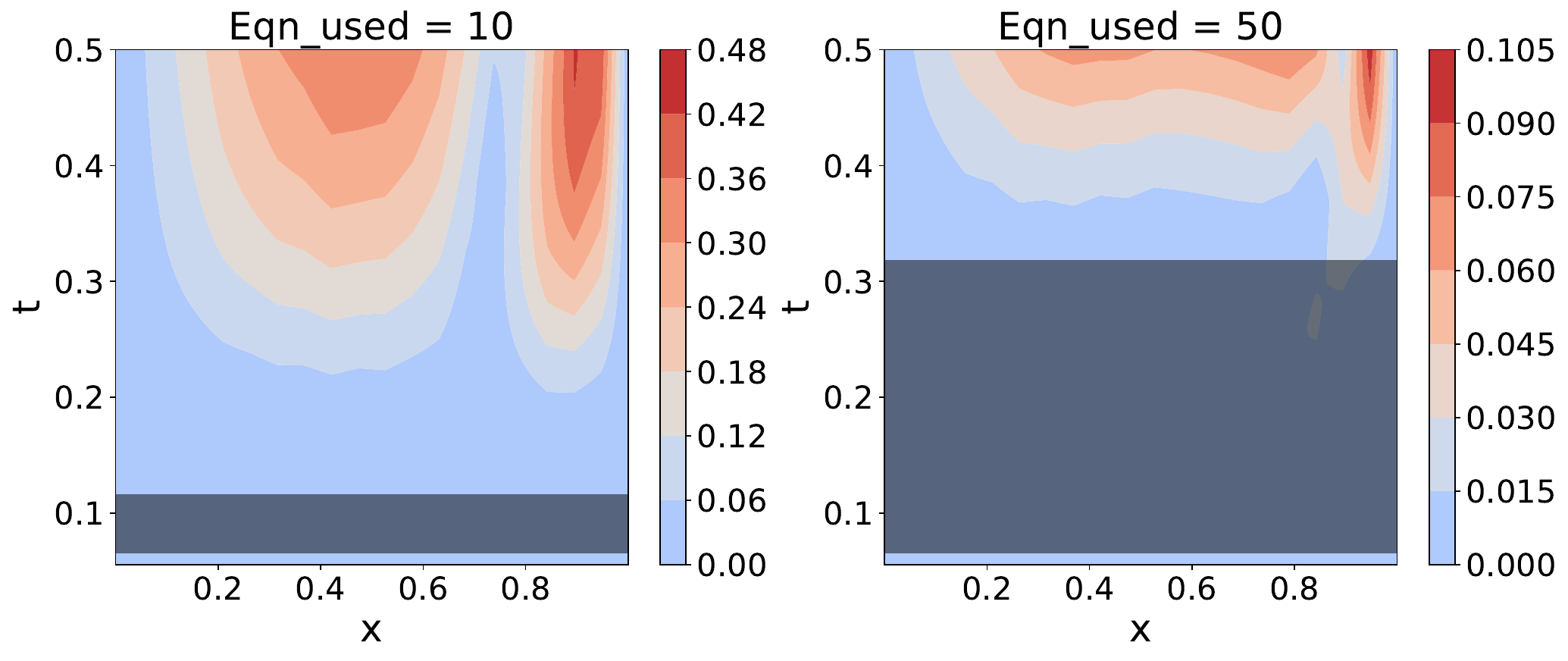}}
\end{subfigure}
\medskip
\caption{Error associated with prediction capabilities when equations solved in the internal spatio-temporal domain and numerical initial point data (=1) used.}
\label{fig:Prediction_Burgers}
\end{figure}

Now, we assess the applicability of the discretized equation-based ANN and LSTM network for prediction at future time steps while the numerical data points and governing equations-based residual minimization are considered only within a few initial time-instants. For the ANN network, we minimize the discretized equation-based residual from $2^{nd}$ time step to $12^{nd}$ time step and from $2^{nd}$ time step to $52^{nd}$ time step with a total of $10$ and $50$ time steps respectively. The initial condition $u = u_{int}$ is the only data-driven constraint considered in the ANN-DisPINN network. As shown in \autoref{fig:pred_burgers1}, when the discretized governing equation \autoref{eq:12} are minimized in the initial $10$ points, the maximum error that occurred in the spatio-temporal domain is $0.35$, whereas in the case of the $50$ initial points, the maximum error value is $0.175$. However, in the case of the LSTM-DisPINN network, we have considered the first $10$ time instants as the sequence length. Therefore, the discretized governing equation-based residuals are minimized from the $12^{th}$ to $22^{nd}$ and from $12^{th}$ to $62^{nd}$ time-step as shown in \autoref{fig:pred_burgers2}. The $u$ value at the $10^{th}$ timesteps is considered as the data-driven constraint in the LSTM-DisPINN network. With $10$ and $50$ time steps, the maximum absolute errors are $0.48$ and $0.105$. Although with LSTM-DisPINN the network prediction accuracy has been improved over the ANN-DisPINN when $50$ time steps are considered, the prediction capability is poorer than the reconstruction capability demonstrated previously where the discretized equation-based residuals are minimized in the entire temporal domain. The prediction capability in terms of absolute error with respect to the benchmark high-fidelity data is demonstrated in \autoref {Table 3}.

\begin{table}[h!] 
\centering
 \begin{tabular}{||c | c | c||} 
 \hline
 Eqns Solved & ANN-DisPINN & LSTM-DisPINN \\ [0.5ex] 
 \hline\hline
 10 & 0.35 & 0.175 \\
 50 & 0.48 & 0.10 \\ [1ex] 
 \hline
 \end{tabular}
 \caption{\label{Table 3} Prediction error associated with the ANN-DisPINN and LSTM-DisPINN for the Burgers' Equation.}
\end{table}

\subsection{Reduced Order Burgers' Equation} \label{subsec:Burgers_reduced}
    
\begin{figure}[ht]
\centering
{\label{fig:Reconserror}\includegraphics[width=1.0\linewidth]{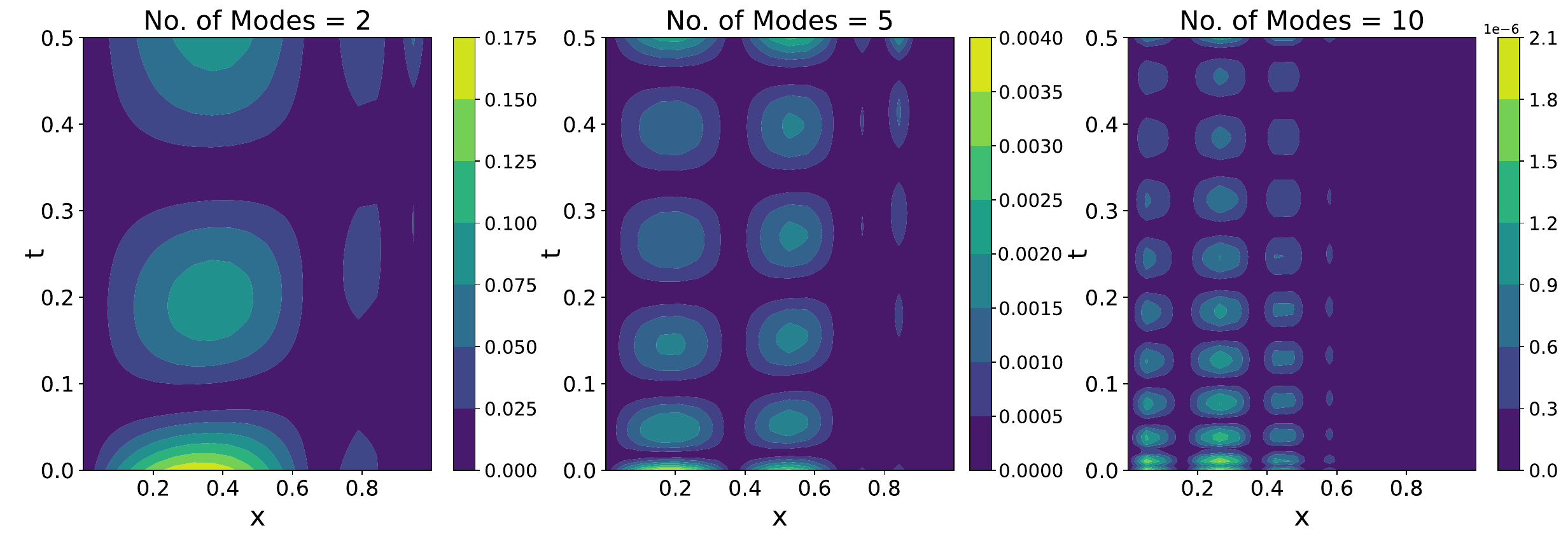}}
\caption{Absolute error associated with the reconstruction with different numbers of POD basis used.}
\label{fig:ReconsError}
\end{figure}

\begin{figure}[ht]
\centering
\begin{subfigure}[b]{0.75\textwidth}
\centering
\subfloat[Reduced  Order ANN-DisPINN]{\label{fig:ReducedANN}\includegraphics[width=1.0\linewidth]{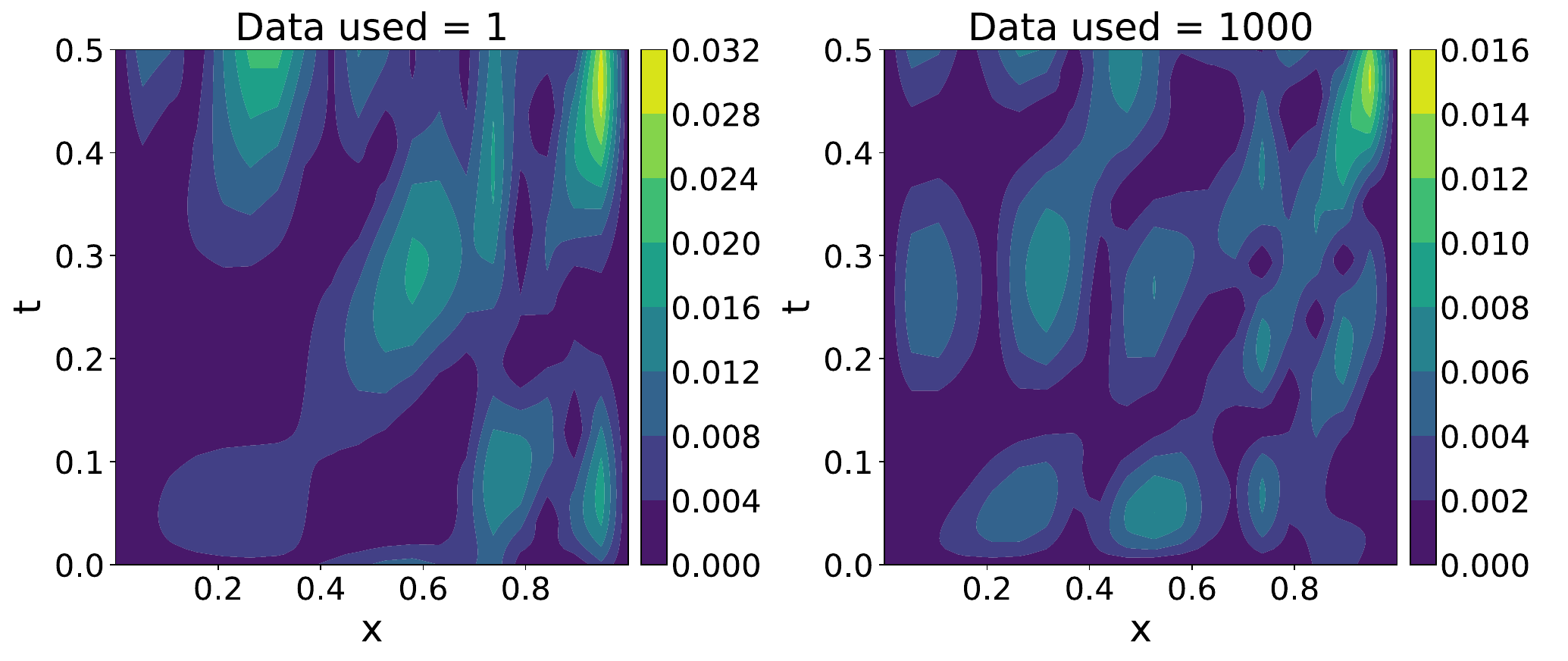}}
\end{subfigure}
\begin{subfigure}[b]{0.75\textwidth}
\centering
\subfloat[Reduced Order 
LSTM-DisPINN]{\label{fig:ReducedLSTM}\includegraphics[width=1.0\linewidth]{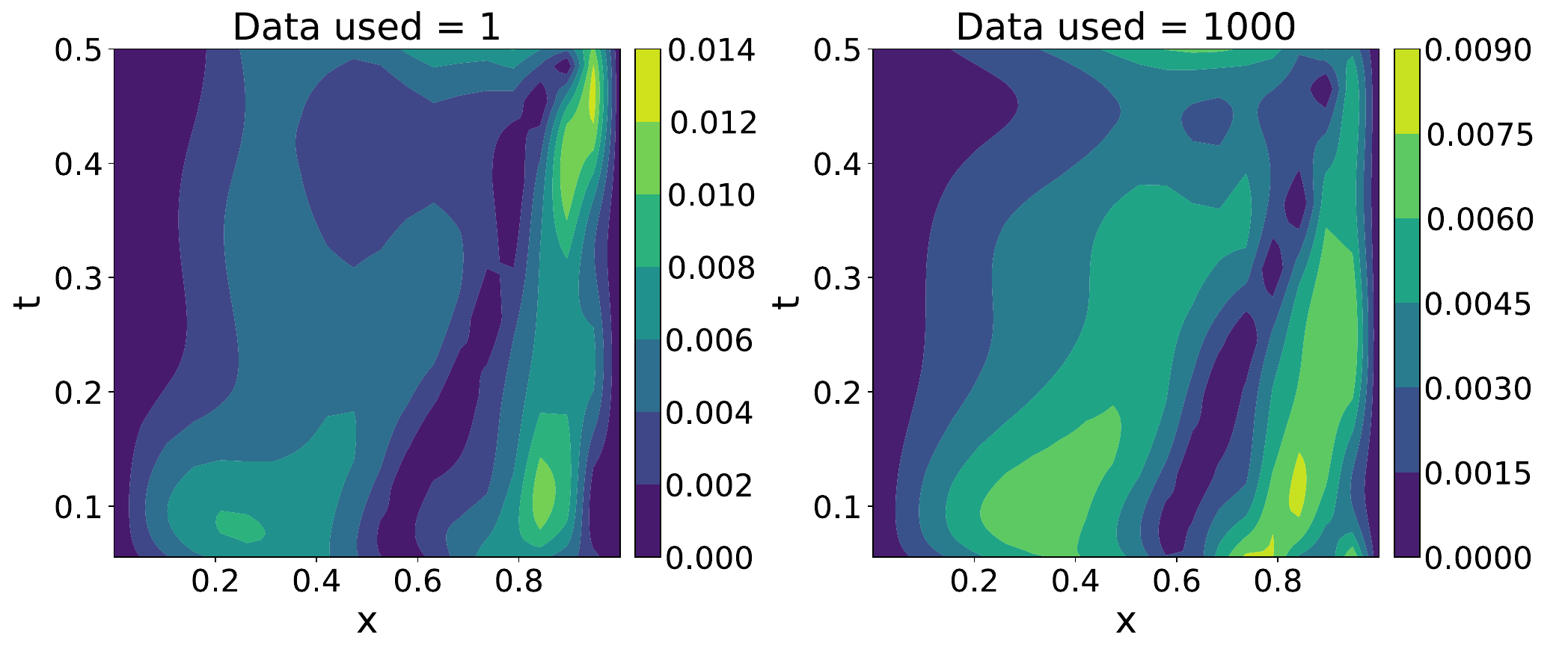}}
\end{subfigure}
\medskip
\caption{Error associated with the reduced order ANN-DisPINN and LSTM-DisPINN with different data points used.}
\label{fig:Reduced}
\end{figure}

In this current section, we first carry out the POD-Galerkin projection of the
Burgers' equation in discretized form, followed by the inclusion of discretized reduced order
equation in the PINN network. $10$ POD modes and $10$ DEIM control points are
considered in the current work. Similar to the previous \autoref{subsec:Burgers_Full}, the entire
temporal domain, which is $100$ time instants, is considered as the input to the neural network. Conversely, in the output, instead of the full discretized equation \autoref{eq:12}, only the reduced order \autoref{eq:13} are considered. Non-linearity arising from the convection term is handled using hyper-reduction as shown in \autoref{eq:17}. Therefore, computational training time savings are achieved. First, the solution vectors at different time instants are considered. Then SVD is applied to the set of solution vectors to compute the POD modes. The reconstruction error is defined as the absolute error between the high-fidelity results and the reconstructed field from the POD coefficients. The reconstruction error is computed with a different number of modes such as $2$, $5$ and $10$ as shown in \autoref{fig:ReconsError}. The maximum error associated with $2$, $5$ and $10$ number of modes are $0.175$, $0.004$ and $10^{(-6)}$ respectively. In the context of a discretized reduced order system-based PINN framework, the reconstruction error first needs to be minimized to expect an overall minimization of error. Hence $10$ modes are considered for the PINN prediction. \autoref{fig:ReducedANN} and \autoref{fig:ReducedLSTM} show the prediction of the PINN network using discretized reduced order governing equation and employing ANN and LSTM as a neural network, respectively. As previously noticed, the LSTM-PINN network outperforms the ANN-PINN network since the maximum error associated with the LSTM-DisPINN and ANN-DisPINN are $0.014$ and $0.032$, respectively, when sparse datasets such as data-point of 1 are considered, and discretized-equation-based residual is minimized in the entire computational domain. However, with the increment of the data points (i.e. $1000$ for ANN-DisPINN and 990 for LSTM-DisPINN), the maximum absolute error drops to $0.016$ and $0.0090$, respectively.\autoref{Table 4} shows the wall clock time required for the training using the full order and reduced order ANN-DisPINN and LSTM-DisPINN while using the $10$ and $5$ DEIM control points and POD modes respectively. It is evident from \autoref{Table 4}, that when the number of modes is reduced to 10 and 5 which is half and a quarter of the total grid size, the computational time is not linearly reduced. Parametric variation and temporal prediction in a future time step will be considered as a future endeavour by the authors to investigate the full potential of the reduced-order discretized physics-based network. Although we intend to reduce the usage of numerical data in ANN-DisPINN and LSTM-DisPINN to demonstrate the potential use of PINN, in the projection phase (POD-Galerkin projection) we have considered the full spatio-temporal domain for the computation of the POD modes and DEIM control points. The authors also aim to reduce the requirements of numerical data in the projection phase using discretized physics-based convolution autoencoder \cite{romor2023non} as a future extension of the current work.     

\begin{table}[h!] 
\centering
 \begin{tabular}{||c | c | c ||} 
 \hline
 No. of Modes and DEIM pts. & ANN-DisPINN & ANN-DisPINN Reduced \\ [1ex] 
 \hline\hline
 10 & 1111.65 $s$ & 937.60 $s$   \\ 
 5 & 1111.65 $s$ & 561.99 $s$ \\ [1ex] 
 \hline
 \hspace{1cm} 
 No. of Modes and DEIM pts.  & LSTM-DisPINN & LSTM-DisPINN Reduced \\ [0.5ex] 
 \hline\hline
 10  & 993.40 $s$ & 852.03 $s$  \\ 
 5 & 993.40 $s$ & 521.52 $s$ \\ [1ex] 
 \hline
 \end{tabular}
 \caption{\label{Table 4} Time taken for different full order (ANN-DisPINN and LSTM-DisPINN) and Reduced Order (ANN-DisPINN and LSTM-DisPINN) training.}
\end{table}

\subsection{Coupling with Detached External Solver} \label{subsec:new_loss}
\begin{figure}[ht]
\centering
\begin{subfigure}[b]{1\textwidth}
\centering
\subfloat[ANN-DisPINN]{\label{fig:jANNh}\includegraphics[width=.50\linewidth]{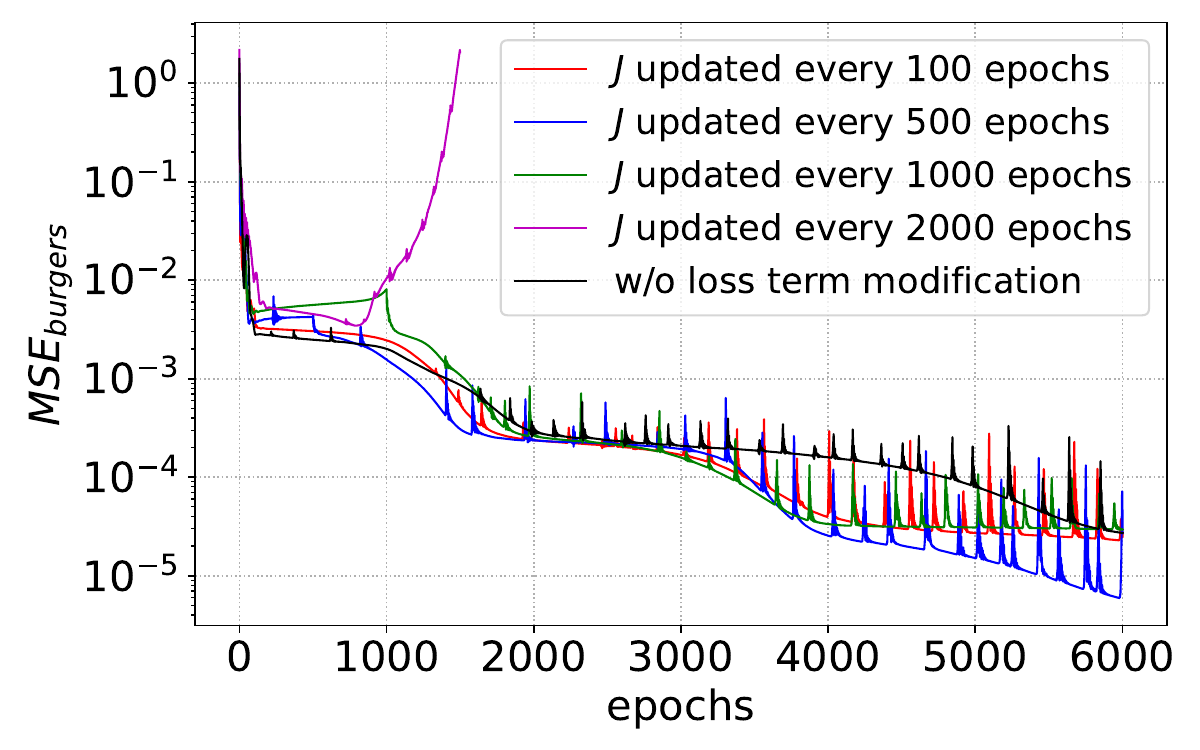}}
\subfloat[LSTM-DisPINN]{\label{fig:jLSTMh}\includegraphics[width=.50\linewidth]{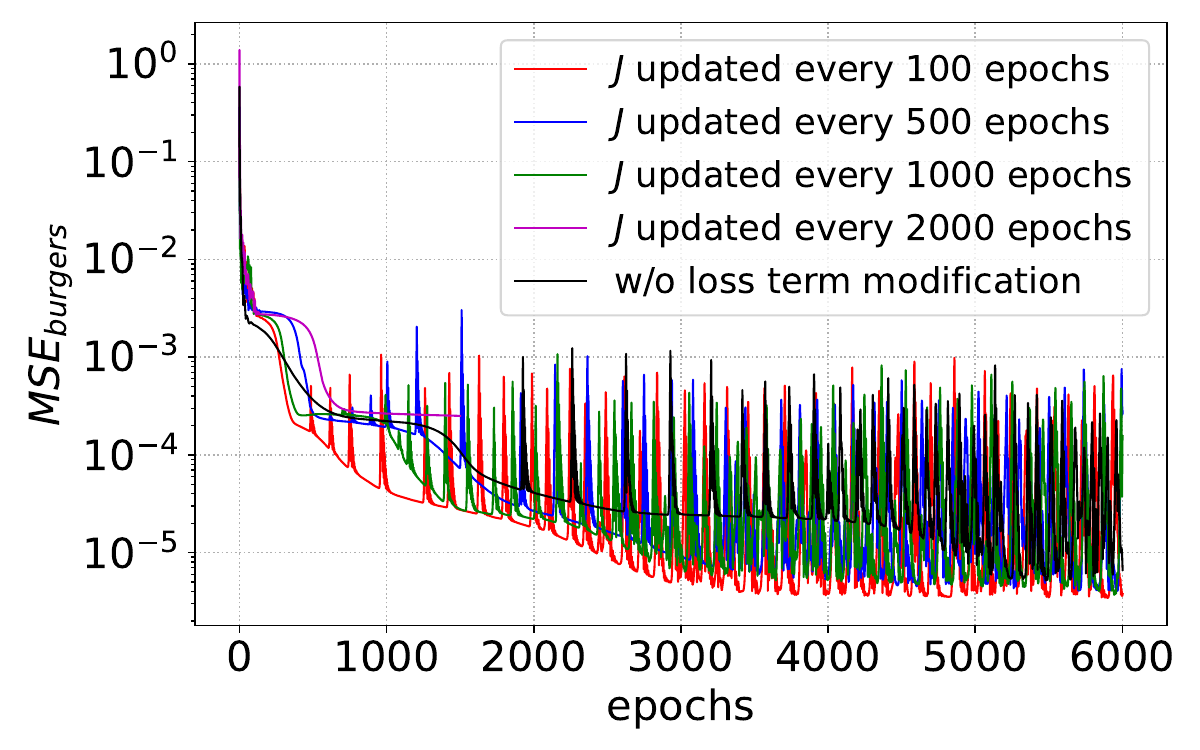}}
\end{subfigure}
\caption{mse associated with the Physics-based loss term with a number of epochs and different jacobian, $J$ update rate.}
\label{fig:Jupdate}
\end{figure}

\begin{figure}[ht]
\centering
{\includegraphics[width=0.85\linewidth]{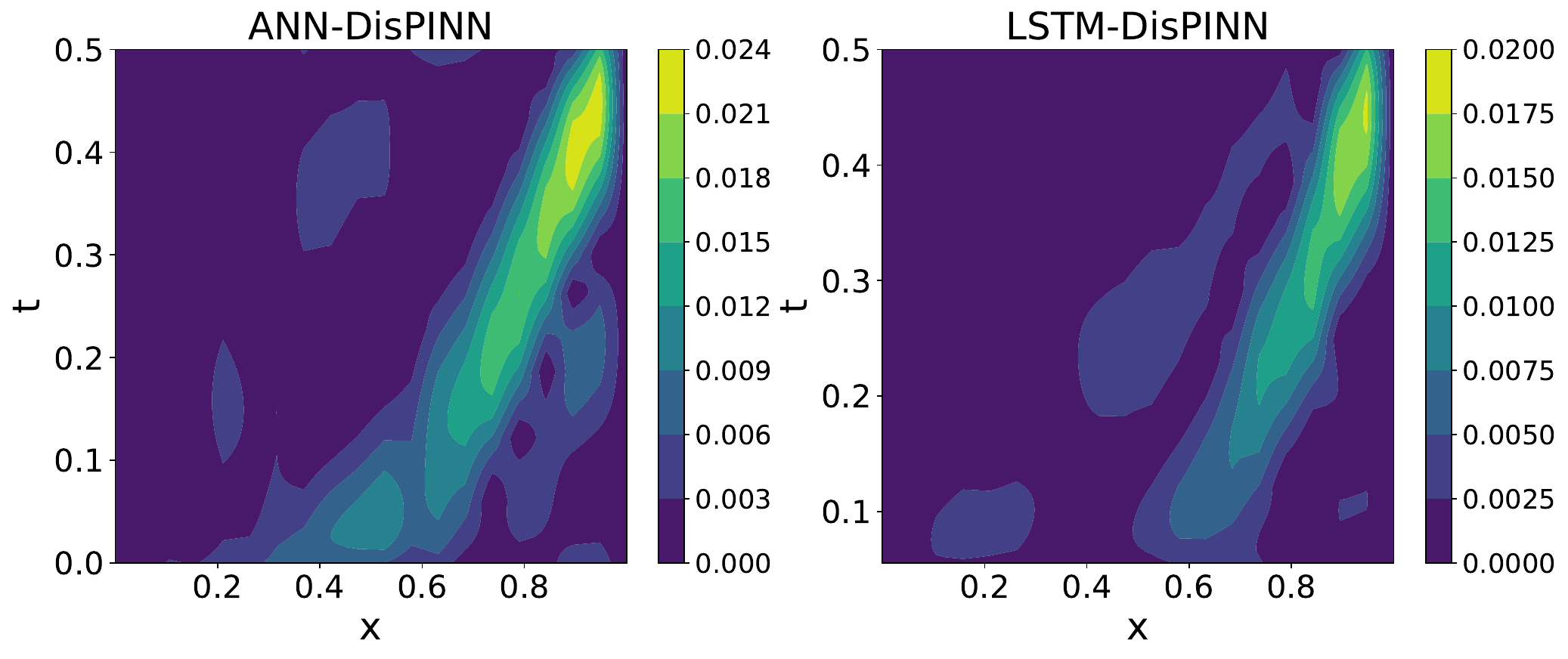}}
\caption{Comparison of the error associated with the ANN-DisPINN and LSTM-DisPINN after new loss term ($\text{L}_{Dis}$) modification (Data used = 1).}
\label{fig:newloss}
\end{figure}

In this subsection, we will discuss the inclusion of the discretized form of the governing equation in the computational graph of the network using the steps mentioned in \autoref{subsec:new_loss} and \autoref{alg:algo}, when the external solver is completely detached from the PINN environment. The external Burgers' equation although implemented in the PyTorch platform, we ensured while passing the residual and the Jacobian, $J$ to the PINN solver, they are detached from the computational graph. This is to ensure, potential coupling with any external solver. \autoref{fig:Jupdate} demonstrates that for ANN-DisPINN, if $J$ is updated at $100$, $500$, $1000$ epochs interval, the governing equation-based residual $MSE_\text{burgers}$ is minimized ensuring the accuracy of the solution at an insignificant increment of the training time. The residual minimization is similar to the case mentioned in \autoref{subsec:Burgers_Full}, where the functional form of the discretized equation is not detached from the computational graph (indicated by black solid lines and termed as w/o loss term modification). If the Jacobian $J$ is updated at an interval of $2000$, in the case of the ANN-DisPINN, the loss term starts increasing after approximately $800$ time-steps. Although the Jacobian update step is expensive, since it is updated only at a few epoch intervals, the total training time is not affected significantly while maintaining the solution accuracy as shown in \autoref{fig:newloss}. Furthermore, in the case of the LSTM-DisPINN, the physics-based loss term is minimized even with $2000$ epoch interval of the $J$ update. The absolute error associated with the Burgers' equation is shown using our proposed new loss term with a Jacobian $J$ update of $500$ and $100$ epoch interval for the ANN-DisPINN and LSTM-DisPINN respectively and the prediction is in good agreement with the benchmark numerical results. The hyper-parameter for the ANN-DisPINN and LSTM-DisPINN is kept similar as mentioned in the \autoref{subsec:Burgers_Full}. Only the initial condition which is $\text{time} \text{step} =0$ for ANN-DisPINN and $\text{time} \text{step} = 10$ for LSTM-DisPINN is considered for the data-driven loss term in \autoref{alg:algo}.

%% file: Conclusion.tex
\section{Conclusion}\label{sec:Conclusion}

In this current work, we propose a novel scheme introducing the discretized governing equations followed by a POD-Galerkin projection of the discretized governing equation and using the numerical residuals in the loss term of the ANN and LSTM network. The proposed approach shows its potential over the conventional Automatic Differentiation-based physics-informed neural network and also lays the foundation for applications in coupling with external forward solvers. The LSTM-PINN network, given an architecture-dependent on the time history of input, possesses several difficulties in terms of the conventional PINN framework. Still, in the current framework, any type of neural network including the LSTM network can be put easily. The additional projection of discretized governing equations onto reduced space improves the training time without compromising the prediction accuracy. Our recent work also demonstrates that if the training dataset is large, the pure data-driven network outperforms the physics-constrained ROM. The introduction of the physics-based loss term appears to deteriorate the performance by increasing the magnitude of the residual. We also propose an efficient algorithm for coupling forward solvers with PINN where it is difficult to access the discretized form of the equation to include them in the computational graph without a significant increase in training time. In future, open-source forward CFD platforms or software dedicated to handling projection-based reduced order schemes will be coupled with the physics-informed neural network framework for complex geometry and large-scale problems.  

%% file: Disclosure_Statement.tex
\section*{\bf{Data Availability}}\label{sec:DA}
The data sets generated during and/or analyzed during the current study are available from the corresponding authors upon reasonable request.

\section*{\bf{Disclosure Statement}}\label{sec:DS}
The authors report no potential conflict of interest.

\section*{Acknowledgements}
The authors gratefully acknowledge the financial support under the scope of the COMET program within the K2 Center “Integrated Computational Material, Process and Product Engineering (IC-MPPE)” (Project No 886385). This program is supported by the Austrian Federal Ministries for Climate Action, Environment, Energy, Mobility, Innovation and Technology (BMK) and for Labour and Economy (BMAW), represented by the Austrian Research Promotion Agency (FFG), and the federal states of Styria, Upper Austria, and Tyrol. The authors also acknowledge financial support from MUR PRIN project $\text{NA\_FROM\_PDEs}$, INDAM GNCS and MUR PRIN project ROMEU.